\documentclass[11pt]{article}
\usepackage{amsmath,amstext,amssymb,amsthm, amsfonts}
\usepackage[title]{appendix}
\usepackage{color}
\newtheorem{theorem}{Theorem}[section]
\newtheorem{lemma}[theorem]{Lemma}

\newtheorem{definition}[theorem]{Definition}

\newtheorem{example}[theorem]{Example}

\newtheorem{corollary}[theorem]{Corollary}

\theoremstyle{remark}
\newtheorem{remark}[theorem]{Remark}

\numberwithin{equation}{section}

\def\Halmos{\mbox{\quad$\square$}}
\def\proof#1{\noindent {\it #1}}
\def\endproof{\vspace{0pt}}

\newcommand{\slim} {\mathop{\rm lim\,sup}}
\newcommand{\ilim} {\mathop{\rm lim\,inf}}
\def\R{\mathbf{R}}
\def\Q{\mathbb{Q}}
\def\olR{{\overline{\R}}}

\def\S{\mathbb{S}}				
\def\SS{\mathbb{S}} 		
\def\SSig{\Sigma} 		
\def\M{\mathcal{M}} 		
\def\P{\mathbf{P}} 		
\def\B{\mathcal{B}}
\def\ee{\varepsilon}

\def\h{\mathbf{I}}
\def\N{\mathbb{N}}			

\def\X{\mathbb{X}}

\def\A{\mathbb{A}}
\def\K{\mathbb{K}}
\def\U{\mathbb{U}}
\def\a{\alpha}				
\def\K{\mathbb{K}}			
\def\PS{\Pi}					
\def\PR{P}					

\begin{document}

\title{Fatou's Lemma in Its Classic Form and Lebesgue's Convergence Theorems for Varying Measures with Applications to MDPs}

\date{}

\author{Eugene~A.~Feinberg\footnote{Department of Applied Mathematics and
Statistics,
 Stony Brook University,
Stony Brook, NY 11794-3600, USA, eugene.feinberg@stonybrook.edu},\
Pavlo~O.~Kasyanov\footnote{Institute for Applied System Analysis,
National Technical University of Ukraine ``Igor Sikorsky Kyiv Polytechnic
Institute'', Peremogy ave., 37, build, 35, 03056, Kyiv, Ukraine,\
kasyanov@i.ua.},\ and Yan~Liang\footnote{Rotman School of Management, University of Toronto,
105 St. George Street, Toronto, ON M5S 3E6, Canada, yan.liang@rotman.utoronto.ca}}

\maketitle

\begin{abstract}
{
The classic Fatou lemma states that the lower limit of a sequence of integrals of functions is greater than or equal to the integral of the lower limit.  It is known that Fatou's lemma for a sequence of weakly converging measures states a weaker inequality because the integral of the lower limit is replaced with the integral of the lower limit in two parameters, where the second parameter is the argument of the functions.  This paper provides sufficient conditions when Fatou's lemma holds in its classic form for a sequence of weakly converging measures.  The functions can take both positive and negative values.  The paper also provides similar results for sequences of setwise converging measures.  It also provides Lebesgue's and monotone convergence theorems for sequences of weakly and setwise converging measures.  The obtained results are used to prove broad sufficient conditions for the validity of optimality equations for average-cost Markov decision processes.
}
\end{abstract}

\section{Introduction}
\label{sec:intr}

For a sequence of nonnegative measurable functions $\{f_n\}_{n\in\N^*},$ Fatou's lemma states the inequality
\begin{equation}\label{eq:Fatcl}
 \int_\SS\liminf_{n\to\infty} f_n(s)\mu(ds) \le \liminf_{n\to\infty}\int_\SS f_n(s)\mu(ds).
\end{equation}
Many problems in probability theory and its applications deal with sequences of probabilities or measures converging in some sense rather than with a single probability or measure $\mu.$  Examples of areas of applications include limit theorems \cite{Bil68}, \cite{JSh}, \cite[Chapter III]{Shi96}, continuity properties of stochastic processes \cite{KL}, and stochastic control \cite{FKZ12, FKZ14, POMDP, HLL96}.

If a sequence of measures $\{\mu_n\}_{n\in\N^*}$ converging setwise to a measure $\mu$ is considered instead of a single measure $\mu,$ then equality \eqref{eq:Fatcl} holds with the measure $\mu$ in its right-hand side replaced with the measures $\mu_n$ \cite[p. 231]{Roy68}.  However, for a sequence of measures $\{\mu_n\}_{n\in\N^*}$ converging weakly to a measure $\mu$, the weaker inequality
	\begin{align}\label{eq:lf}
		\int_\S \ilim\limits_{n\to\infty,\, s'\to s} f_n (s')\mu(ds)
		\le \ilim\limits_{n\to\infty}\int_\S f_n (s)\mu_n (ds)
	\end{align}
holds.  Studies of Fatou's lemma for weakly converging probabilities were started by Serfozo~\cite{Ser82} and continued in \cite{FKL18, FKZTVP}. For a sequence of measures converging in total variation, Feinberg et al.~\cite{FKZ16} obtained the uniform Fatou's lemma, which is a more general fact than Fatou's lemma.

This paper describes sufficient conditions ensuring that Fatou's lemma holds in its classic form for a sequence of weakly converging measures.  In other words, we provide sufficient conditions for the validity of inequality \eqref{eq:ilimn}, which is inequality \eqref{eq:lf} with its left-hand side replaced with the left-hand side of \eqref{eq:Fatcl}.  We consider the sequence of functions that can take both positive and negative values.  In addition to the results for weakly converging measures, we provide parallel results for setwise converging measures.  We also investigate the validity of Lebesgue's and monotone convergence theorems for sequences of weakly and setwise converging measures.  The results are applied to Markov decision processes (MDPs) with long-term average costs per unit time, for which we provide general conditions for the validity of optimality equations.

Section~\ref{sec:intro:fatou} describes the three types of convergence of measures: weak convergence, setwise convergence, and convergence in total variation, and it provides the known formulations of Fatou's lemmas for these types of convergence modes.  Section~\ref{sec:auxil} describes conditions under which the double lower limit of a sequence of functions in the left-hand side of \eqref{eq:lf} is equal to the standard lower limit.  Section~\ref{sec:Fatou} describes sufficient conditions for the validity of Fatou's lemma in its classic form for a sequence of weakly converging measures.  This section also provides results for sequences of  measures converging setwise.  Sections~\ref{sec:main} and \ref{sec:swtv} describe Lebesgue's and monotone convergence theorems for  weakly and setwise converging measures.  Section~\ref{sec:appl} deals with applications. 

\section{Known Formulations of Fatou's Lemmas for Varying Measures}
\label{sec:intro:fatou}

Let $(\SS,\SSig)$ be a measurable space, $\M (\SS)$ be the {\it family of all 
finite measures on} $(\SS,\SSig),$ and $\P(\SS)$ be the {\it family of all probability measures} on
$(\SS,\SSig).$ When $\SS$ is a metric space, we always consider
$\SSig := \B(\SS),$ where $\B(\SS)$ is the {\it Borel $\sigma$-field} on $\SS.$
Let $\R$ be the real line, $\olR := [-\infty,+\infty],$ and $\N^* :=\{1,2,\ldots\}.$
We denote by $\h \{ A \}$  the \textit{indicator} of the event $A.$

Throughout this paper, we deal with integrals of functions
that can take both positive and negative values.  An \textit{integral} $\int_\SS f (s)\mu(ds)$
of a measurable $\overline{\R}$-valued function $f$ on $\SS$ with
respect to a measure $\mu$ is \textit{defined} if
\begin{align}\label{e:condint}
\min\{ \int_\SS f^+(s)\mu(ds),
\int_\SS f^-(s) \mu(ds)\}< +\infty,
\end{align}
where $f^+(s)=\max\{f(s),0\},$ $f^-(s)=-\min\{f(s),0\},$ $s\in \SS.$  If
\eqref{e:condint} holds, then the integral  is defined as $\int_\SS f
(s)\mu(ds)=\int_\SS f^+ (s)\mu(ds)- \int_\SS f^- (s)\mu(ds).$  All the
integrals in the assumptions of the following lemmas, theorems, and
corollaries are assumed to be defined.
For $\mu \in \M (\SS)$
consider the vector space $L^1(\SS; \mu)$ of all measurable functions $f : \SS \mapsto \olR,$
whose absolute values have finite integrals, that is, $\int_\SS |f(s)|\mu(ds) < +\infty.$

We recall the definitions of the following three types of convergence of measures:  weak convergence, setwise convergence, and convergence in total variation.

{
\begin{definition}[Weak convergence]\label{def:wc}
	A sequence of measures $\{ \mu_n \}_{n\in\N^*}$ on a metric space $\S$
	\textit{converges weakly} to a finite measure $\mu$ on $\S$ if, for each bounded continuous function $f$ on $\S,$
		 \begin{align}
	\int_\S f(s) \mu_n (ds) \to \int_\S f(s) \mu (ds)  \quad \text{ as } n\to\infty.
	\label{eq:def:setwise}
\end{align}
\end{definition}


\begin{remark}\label{rem:fin meas weak}
{\rm  Definition~\ref{def:wc} implies that
 $\mu_n(\S)\to \mu(\S)\in \R$ as $n\to\infty.$  Therefore, if $\{ \mu_n \}_{n\in\N^*}$ converges weakly to $\mu\in  \M (\S),$ then
there exists $N\in\N^*$ such that $\{\mu_n\}_{n=N,N+1,\ldots} \subset \M (\S).$ 
}
\end{remark}


\begin{definition}[Setwise convergence]\label{def:sw}
	A sequence of measures $\{ \mu_n \}_{n\in\N^*}$ on a measurable space $(\SS,\SSig)$
	\textit{converges setwise} to a measure $\mu$ on $(\SS,\SSig)$
	if for each $C\in\SSig$
	\begin{align*}
		 \mu_n (C) \to \mu (C)  \quad \text{ as } n\to\infty.
	\end{align*}
\end{definition}

\begin{definition}[Convergence in total variation]\label{def:tv}
	A sequence of finite measures $\{ \mu_n \}_{n\in\N^*}$ on a measurable space $(\SS,\SSig)$
	\textit{converges in total variation} to a measure $\mu$ on $(\SS,\SSig)$ if
	\begin{align*}
		\sup\Bigg\{ \Big| \int_\SS f(s) \mu_n (ds) - \int_\SS f(s) \mu (ds) \Big|\, :\, f:\SS\mapsto [-1,1] \text{ is measurable}  \Bigg\} \to 0 \quad \text{ as } n\to\infty.
	\end{align*}
\end{definition}    }

\begin{remark}\label{rem:tv setwise weak}
{\rm
As follows from Definitions~\ref{def:wc}, \ref{def:sw}, and \ref{def:tv},  if a sequence of finite measures $\{ \mu_n \}_{n\in\N^*}$ on a measurable space $(\SS,\SSig)$
converges in total variation to a measure $\mu$ on $(\SS,\SSig),$ then $\{ \mu_n \}_{n\in\N^*}$ converges setwise to $\mu$ as $n\to\infty$ and the measure  $\mu$ is finite.  The latter follows from $|\mu_n(\SS)-\mu(\SS)|<+\infty$ when $n\ge N$ for some $N\in\N^*.$ Furthermore, if a sequence of measures $\{ \mu_n \}_{n\in\N^*}$ on a metric space $\S$ converges setwise to a finite measure $\mu$ on $\S,$ then this sequence converges weakly to $\mu$ as $n\to\infty.$
}
\end{remark}


Recall the following definitions  of the  uniform and asymptotic uniform integrability of sequences of functions.

\begin{definition}
\label{defi:unif integr}
The sequence $\{f_n\}_{n\in\N^*}$  of measurable $\overline{\R}$-valued functions is called
\begin{itemize}
\item {\it uniform integrable (u.i.) with respect to (w.r.t.) a sequence of measures $\{\mu_n\}_{n\in\N^*}$} if
\begin{align}
		\lim_{K\to+\infty} \sup_{n\in\N^*} \int_\SS |f_n (s)|
		\h \{ s\in\SS : |f_n (s)| \geq K\} \mu_n (ds) = 0 ;
		\label{eq:tv:ui}
	\end{align}
\item {\it asymptotically uniform integrable (a.u.i.) w.r.t.  a sequence of measures  $\{\mu_n\}_{n\in\N^*}$}  if
\begin{align}
		\lim_{K\to+\infty} \slim_{n\to\infty} \int_\SS |f_n (s)|
		\h \{ s\in\SS : |f_n (s)| \geq K\} \mu_n (ds) = 0 .
		\label{eq:tv:ui:lim}
	\end{align}
\end{itemize}
\end{definition}

If $\mu_n=\mu\in\M(\SS)$ for each $n\in\N^*,$ then an (a.)u.i.~w.r.t. $\{\mu_n\}_{n\in\N^*}$
sequence $\{f_n\}_{n\in\N^*}$ is called (\textit{a.})\textit{u.i.} For  $\mu\in\M(\SS)$ a sequence $\{f_n\}_{n\in\N^*}$ of functions from $L^1 (\SS;\mu)$
is u.i. if and only if it is a.u.i.; Kartashov~\cite[p. 180]{Kar}.
For a single finite measure $\mu,$ the definition of an a.u.i.~sequence of functions
(random variables in the case of a probability measure $\mu$) coincides with the corresponding
definition broadly used in the literature; see, e.g., \cite[p. 17]{Vaa98}. Also, for a single fixed
finite measure, the definition of a u.i. sequence of functions is consistent with the classic
definition of a family $\cal H$ of u.i. functions.
We say that a function $f$ is (a.)u.i.~w.r.t.
$\{\mu_n\}_{n\in\N^*}$ if the sequence $\{f,f,\ldots\}$ is  (a.)u.i.~w.r.t.
$\{\mu_n\}_{n\in\N^*}.$ A function $f$ is  u.i.~w.r.t. a family  $ \mathcal{N} $
of measures if
\begin{align*}
	\lim_{K\to+\infty} \sup_{\mu\in  \mathcal{N}  }  \int_\SS |f (s)|
	\h \{ s\in\SS : |f (s)| \geq K\} \mu (ds) = 0.
\end{align*}


\begin{theorem}[Equivalence of u.i. and a.u.i.; {Feinberg et al.~\cite[Theorem 2.2]{FKL18}}]\label{thm:uiCondEqui}
Let $(\SS, \SSig)$ be a measurable space, $\{\mu_n\}_{n\in\N^*} \subset \M (\SS),$ and
$\{ f_n \}_{n\in\N^*}$ be a sequence of measurable
$\olR$-valued functions on $\SS.$ Then there exists $N\in\N^*$ such that $\{f_n\}_{n=N,N+1,\ldots}$ is u.i. w.r.t. $\{\mu_{n}\}_{n=N,N+1,\ldots}$ iff $\{f_n\}_{n\in\N^*}$ is a.u.i. w.r.t. $\{\mu_n\}_{n\in\N^*}.$
\end{theorem}


Fatou's lemma (FL) for weakly converging  probabilities was introduced in Serfozo~\cite{Ser82} and generalized in \cite{FKL18, FKZTVP}.

\begin{theorem}[{FL for weakly converging measures; \cite[Theorem 2.4 and Corollary 2.7]{FKL18}}]\label{thm:Fatou:w}
	Let $\S$ be a metric space, $\{\mu_n\}_{n\in\N^*}$ be a sequence of measures on $\S$ that converges weakly to $\mu\in \M (\S),$ and
	$\{ f_n \}_{n\in\N^*}$ be a sequence of measurable $\olR$-valued functions on $\S.$
	Assume that one of the following two conditions holds:
	\begin{itemize}
		\item[{\rm (i)}] $\{ f^-_n \}_{n\in\N^*}$ is a.u.i. w.r.t. $\{ \mu_n\}_{n\in\N^*};$
		\item[{\rm (ii)}] there exists a sequence of measurable real-valued functions $\{g_n\}_{n\in\N^*}$
	on $\S$ such that $ f_n (s)\ge g_n(s)$ for all $n\in\N^*$ and $s\in \S,$ and
	\begin{align}\label{eq:sw2aaa:var}
		-\infty<\int_\S\slim_{n\to\infty,s^\prime\to s} g_n(s^\prime)\mu(ds)
		\le\ilim_{n\to\infty}\int_\S g_n(s)\mu_n(ds) .
	\end{align}
	\end{itemize}
	  Then inequality \eqref{eq:lf} holds. 
\end{theorem}


Recall that FL for setwise converging measures is stated in Royden~\cite[p. 231]{Roy68} for nonnegative functions. FL for setwise converging probabilities is stated in Feinberg~et~al.~\cite[Theorem~4.1]{FKZTVP} for  functions taking positive and negative values.


\begin{theorem}[{FL for setwise converging probabilities; \cite[Theorem~4.1]{FKZTVP}}]\label{thm:Fatou:sw}
Let $(\SS,\SSig)$ be a measurable space, a sequence of measures
$\{\mu_n\}_{n\in\N^*}\subset\P(\S)$ converge setwise to $\mu\in\P(\S),$
and $\{ f_n \}_{n\in\N^*}$ be a sequence of measurable real-valued
functions on $\SS.$ Then the inequality
\begin{align}
	\int_\SS \ilim\limits_{n\to\infty} f_n (s)\mu(ds)\le \ilim\limits_{n\to\infty}\int_\SS  f_n (s)\mu_n (ds)
	\label{eq:ilimn}
\end{align}
holds, if there exists a sequence of measurable real-valued functions $\{g_n\}_{n\in\N^*}$
	on $\SS$ such that $ f_n (s)\ge g_n(s)$ for all $n\in\N^*$ and $s\in \SS,$ and
	\begin{align}\label{eq:sw2aaa:sw:var}
		-\infty<\int_\SS\slim_{n\to\infty} g_n(s)\mu(ds)\le\ilim_{n\to\infty}\int_\SS g_n(s)\mu_n(ds) .
	\end{align}
\end{theorem}


Under the condition that
$\{\mu_n\}_{n\in\N^*} \subset \M (\SS)$ converges in total variation
to $\mu\in\M(\SS),$ Feinberg et al.~\cite[Theorem~2.1]{FKZ16} established uniform FL,
which is a stronger statement than the classic FL. 


\begin{theorem}[{Uniform FL for measures converging in total variation; \cite[Theorem~2.1]{FKZ16}}]\label{thm:UFL:FKZ}
Let $(\SS, \SSig)$ be a measurable space, a sequence of measures
$\{\mu_n\}_{n\in\N^*}$ from $\M (\SS)$ converge
in total variation to a measure $\mu\in\M(\SS),$
$\{ f_n\}_{n\in\N^*}$ be a sequence of measurable $\olR$-valued functions on $\SS,$ and $f$ be a measurable $\olR$-valued function.
Assume that $f\in L^1 (\SS;\mu)$ and $f_n\in L^1 (\SS;\mu_n)$ for each $n\in\N^*.$ Then the inequality
\begin{align}
	\ilim_{n\to\infty} \inf_{C\in\SSig} \Big( \int_C f_n (s)\mu_n (ds) - \int_C f(s) \mu (ds) \Big) \geq 0
	\label{eq:ufl:tv}
\end{align}
holds if and only if the following two statements hold:
\begin{itemize}
	\item[{\rm (i)}] for each $\ee > 0$
	$
		\mu (\{ s\in\SS : f_n (s) \leq f(s) - \ee \}) \to 0 \text{ as } n\to \infty,
		$
and, therefore, there exists a subsequence $\{f_{n_k}\}_{k\in\N^*} \subset \{f_n\}_{n\in\N^*}$
such that
$f(s) \le \ilim_{k\to\infty} f_{n_k} (s)$
for $\mu$-a.e. $s\in\SS;$
\item[{\rm (ii)}] $\{f_n^-\}_{n\in\N^*}$ is a.u.i. w.r.t. $\{\mu_n\}_{n\in\N^*}.$
\end{itemize}
\end{theorem}

\section{Semi-Convergence Conditions for  Sequences of Functions}
\label{sec:auxil}

Let $(\SS,\SSig)$ be a measurable space, $\mu$ be a measure on $(\SS,\SSig),$ $\{f_n\}_{n\in\N^*}$ be a sequence of measurable $\olR$-valued functions, and $f$ be a measurable $\olR$-valued function.
In this section we establish the notions of
lower and upper semi-convergence in measure $\mu$ (see Definition~\ref{def:lscvgm}) for a sequences of functions $\{f_n\}_{n\in\N^*}$ defined on a measurable space $\S$. Then, under the assumption that $\S$ is a metric space, we examine necessary and sufficient conditions for the following equalities (see Theorem~\ref{lm:EC},  Corollary~\ref{cor:limit}, and Example~\ref{ex:equi}):
\begin{align}
\ilim_{n\to \infty,s^\prime\rightarrow s}  f_n  (s^\prime)&=\ilim_{n\to\infty}  f_n (s),
		\label{eq:llim}\\
\lim_{n\to \infty,s^\prime\rightarrow s}  f_n  (s^\prime) &= \lim_{n\to\infty}  f_n (s),
		\label{eq:lim}
	\end{align}
which improve the statements of FL  and Lebesgue's  convergence theorem for weakly converging measures; see Theorem~\ref{thm:Fatou:w:lsec} and Corollary~\ref{cor:leb:weak}.  For example, these equalities are important for approximating average-cost relative value functions for MDPs with weakly continuous transition probabilities by discounted relative value functions; see Section~\ref{sec:appl}.
For this purpose we introduce the notions of lower and upper semi-equicontinuous families of functions; see Definition~\ref{def:lowerequicont}. Finally, we provide  sufficient conditions for lower semi-equicontinuity; see Definition~\ref{defi:unconverg} and  Corollary~\ref{cor:usc-lsec}.


\begin{remark}\label{rem:mainconv}
{\rm
Since
\begin{equation}\label{eq:P1}
\ilim_{n\to \infty,s^\prime\rightarrow s}  f_n  (s^\prime) \le \ilim_{n\to\infty} f_n (s),
\end{equation}
 \eqref{eq:llim} is equivalent to the  inequality
\begin{equation}\label{eq:P2}
\ilim_{n\to\infty} f_n (s)\le \ilim_{n\to \infty,s^\prime\rightarrow s}  f_n  (s^\prime).
\end{equation}
}
\end{remark}


 To provide  sufficient conditions for (\ref{eq:llim}) we introduce the
 definitions of  uniform semi-convergence.

\begin{definition}[Uniform semi-convergence]\label{defi:unconverg}
A sequence of
real-valued functions $\{f_n\}_{n\in\N^*}$ on $\S$
semi-converges uniformly from below to a real-valued function $f$ on $\S$ if for each $\ee > 0$ there exists $N\in\N^*$ such that
	\begin{align}
		f_n (s) > f(s) - \ee
		\label{eq:usc1}
	\end{align}
for each $s\in\S$  and $n =N,N+1,\ldots\ .$
A sequence of real-valued functions $\{f_n\}_{n\in\N^*}$ on $\S$
semi-converges uniformly from above to a real-valued function $f$ on $\S$ if $\{-f_n\}_{n\in\N^*}$
{semi-converges uniformly from below} to  $-f$ on $\S.$
\end{definition}

\begin{remark}\label{rem:semi-conv}
{\rm  A sequence $\{f_n\}_{n\in\N^*}$ \textit{converges uniformly  to} $f$ on $\S$ if and only if it uniformly semi-converges from below and from above.
}
\end{remark}

Let us consider the following definitions of semi-convergence in measure.
\begin{definition}[Semi-convergence in measure]\label{def:lscvgm}
A sequence of measurable $\overline{\R}$-valued  functions $\{f_n\}_{n\in\N^*}$	lower semi-converges to a measurable real-valued  function $f$
	in measure $\mu$
	if  for each $\ee > 0$
	\begin{align*}
		\mu (\{ s\in\SS : f_n (s) \leq f(s) - \ee \}) \to 0 \text{ as } n\to \infty .
		\end{align*}
		A sequence of measurable $\overline{\R}$-valued  functions $\{f_n\}_{n\in\N^*}$	upper semi-converges to a measurable real-valued  function $f$
	in measure $\mu$
	if $\{-f_n\}_{n\in\N^*}$ lower semi-converges to $-f$ in measure $\mu,$ that is,
for each $\ee > 0$
	\begin{align*}
		\mu (\{ s\in\SS : f_n (s) \geq f(s) + \ee \}) \to 0 \text{ as } n\to \infty .
			\end{align*}
\end{definition}


\begin{remark}\label{rem:conv in meas}
{\rm   A sequence of measurable $\overline{\R}$-valued  functions $\{f_n\}_{n\in\N^*}$ converges to a measurable real-valued  function $f$ \textit{in measure} $\mu,$ that is, for each $\ee > 0$
\begin{align*}
		\mu (\{ s\in\SS : |f_n (s) - f(s)| \ge \ee \}) \to 0 \text{ as } n\to \infty ,
			\end{align*}
if and only if this sequence of functions both lower and upper semi-converges to $f$ in measure $\mu$.  }
\end{remark}


\begin{remark}\label{rem:pointwise semiconv}
{\rm
If 
$f(s) \le \ilim_{n\to\infty} f_n (s),$ $f(s) \ge \slim_{n\to\infty} f_n (s),$
 or $f(s) = \lim_{n\to \infty} f_n (s)$    for    $\mu$-a.e.  $s\in\SS,$
then $\{f_n \}_{n\in\N^*}$ lower semi-converges, upper semi-converges, or converges respectively
to $f$ in measure $\mu.$
Conversely, Feinberg~et~al.~\cite[Lemma 3.1]{FKZ16} implies that if $\{f_n \}_{n\in\N^*}$ lower semi-converges, upper semi-converges, or converges to $f$ in measure $\mu,$ then there exists a subsequence $\{f_{n_k}\}_{k\in\N^*} \subset \{f_n\}_{n\in\N^*}$
such that 
$f(s) \le \ilim_{k\to\infty} f_{n_k} (s),$ $f(s) \ge \slim_{k\to\infty} f_{n_k} (s),$
 or $f(s) = \lim_{k\to \infty} f_{n_k} (s)$  respectively  for    $\mu$-a.e.  $s\in\SS.$

}
\end{remark}


Now let $\S$ be a metric space and
$B_\delta(s)$ be the open ball in $\S$ of radius $\delta>0$ centered at $s\in \S.$ 
We consider the notions of lower and upper semi-equicontinuity for a sequence of functions.

\begin{definition}[Semi-equicontinuity]\label{def:lowerequicont}
	A sequence $\{f_n\}_{n\in\N^*}$ of real-valued functions on a metric space $\S$
	is called \textit{lower semi-equicontinuous at the point} $s\in \S$ if for each $\ee > 0$
	there exists $\delta>0$ such that
	\begin{equation*}
		 f_n(s^\prime) >  f_n(s) - \ee \qquad \text{for all}\  s^\prime\in B_\delta(s) \text{ and for all } n\in\N^*.
	\end{equation*}
	The sequence $\{f_n\}_{n\in\N^*}$ is called \textit{lower semi-equicontinuous} (\textit{on $\S$}) if it is lower semi-equicontinuous at
	all $s \in \S.$
	A sequence $\{f_n\}_{n\in\N^*}$ of real-valued functions on a metric space $\S$
	is called \textit{upper semi-equicontinuous at the point} $s\in \S$ (\textit{on $\S$})
if the sequence $\{ -f_n \}_{n\in\N^*}$ is lower semi-equicontinuous at the point $s\in \S$ (on~$\S$).
\end{definition}

Recall the definition of equicontinuity of a sequence of functions;
see e.g. Royden~\cite[p. 177]{Roy68}.


\begin{definition}[Equicontinuity]\label{def:quicont}
A sequence $\{f_n\}_{n\in\N^*}$  of real-valued functions on a metric space $\S$
is called \textit{equicontinuous at the point} $s\in \S$ (\textit{on $\S$}) if this sequence is
both lower and upper semi-equicontinuous at the point $s\in \S$ (on $\S$).
\end{definition}


Theorem~\ref{lm:EC} states necessary and sufficient conditions for equality (\ref{eq:llim}).
This theorem and Corollary~\ref{cor:limit} generalize Feinberg and Liang~\cite[Lemma 3.3]{FLi17}, where equicontinuity was considered.


\begin{lemma}\label{lm:mono:lim}
Let $\{  f_n  \}_{n\in\N^*}$ be a pointwise nondecreasing sequence of lower semi-continuous $\olR$-valued functions on a metric space $\S.$ Then
\begin{align}
\ilim_{n\to \infty,\,s^\prime\rightarrow s}  f_n  (s^\prime)&=\lim_{n\to\infty}  f_n (s), \qquad s\in\S.
		\label{eq:llim:lim}
\end{align}
\end{lemma}


\proof{Proof.} For each $s\in\S,$
\begin{align*}
	\ilim_{n\to\infty,s^\prime\to s} f_n (s^\prime) = \sup_{n\geq 1}\ilim_{s^\prime\to s}\inf_{k\geq n} f_k (s^\prime) = \sup_{n\geq 1}\ilim_{s^\prime\to s} f_n (s^\prime) = \sup_{n\geq 1} f_n (s)  = \lim_{n\to\infty} f_n (s),
\end{align*}
where the first equality follows from the definition of $\ilim,$ the third one follows from the lower semi-continuity of the function $f_n,$ and the second and the last equalities hold because the sequences $\{f_n \}_{n\in\N^*}$ are pointwise nondecreasing. Hence, \eqref{eq:llim:lim} holds.
\hfill\Halmos\endproof

\begin{theorem}[Necessary and sufficient conditions for (\ref{eq:llim})]\label{lm:EC}
Let $\{  f_n  \}_{n\in\N^*}$ be a sequence of real-valued functions on a
metric space $\S,$ and let $s\in\S.$ Then the following statements hold:
\begin{itemize}
	\item[{\rm (i)}] if the sequence of functions $\{  f_n  \}_{n\in\N^*}$ is lower semi-equicontinuous at $s,$ then each function $f_n,$ $n\in\N^*,$ is lower semi-continuous at $s$ and equality	(\ref{eq:llim}) holds;
		\item[{\rm (ii)}]  if $\{  f_n  \}_{n\in\N^*}$ is the sequence of lower semi-continuous functions  satisfying  \eqref{eq:llim} and $\{ f_n (s) \}_{n\in\N^*}$ is
	a converging sequence, that is,
\begin{equation}\label{eq:eq}
\ilim_{n\to\infty}  f_n (s)=\slim_{n\to\infty}  f_n (s),
\end{equation}
then the sequence $\{  f_n  \}_{n\in\N^*}$ is lower semi-equicontinuous at $s.$
\end{itemize}
\end{theorem}

Example~\ref{ex:equi} demonstrates that Assumption~(\ref{eq:eq}) is essential in Theorem~\ref{lm:EC}(ii). Without this assumption, the remaining conditions of Theorem~\ref{lm:EC}(ii) imply only the existence of a subsequence $\{  f_{n_k}  \}_{k\in\N^*}\subset\{  f_n  \}_{n\in\N^*}$ such that $\{  f_{n_k}  \}_{k\in\N^*}$ is lower semi-equicontinuous at $s.$ This is true because every  subsequence $\{  f_{n_k}  \}_{k\in\N^*}$ satisfying $\lim_{k\to\infty} f_{n_k} (s) = \ilim_{n\to\infty} f_n (s)$ is lower semi-equicontinuous at $s$ in view of Theorem~\ref{lm:EC}(ii) since \eqref{eq:eq} holds for such subsequences.




\begin{example}\label{ex:equi}
{\rm
Consider $\S := [-1,1]$ endowed with the standard Euclidean metric and
	\begin{align*}
		f_n (t) :=
		\begin{cases}
			0, & \text{if } n = 2k-1,  \\
			\max\{ 1 -  n|t|, 0 \}, & \text{if } n = 2k,
		\end{cases}
		\qquad  k\in\N^*, \quad t\in\S .
	\end{align*}
Each function $f_n,$ $n\in\N^*,$ 
is
nonnegative and continuous  on $\S .$
Equality \eqref{eq:llim} holds because
\[ 0\le
\ilim_{n\to \infty,s^\prime\rightarrow 0}  f_n  (s^\prime)
\le\ilim_{n\to\infty}  f_n (0)= f_{2k-1}(0)=0,\qquad k\in\N^*.
\]
Equality \eqref{eq:eq} does not hold because
\[
\slim_{n\to \infty}  f_n  (0)=1>0=\ilim_{n\to \infty}  f_n  (0),
\]
where the first equality holds because $f_{2k}(0)=1$ for each $k\in\N^*,$ and the second equality holds because $f_{2k-1}(0)=0$ for each $k\in\N^*.$
The sequence of functions $\{f_n\}_{n\in\N^*}$ is not lower semi-equicontinuous at $s=0$ because
$f_{2k}(\frac{1}{2k})=0<\frac12=f_{2k}(0)-\frac12$ for each $k\in\N^*.$
	Therefore, the conclusion of Theorem~\ref{lm:EC}(ii) does not hold, and assumption~(\ref{eq:eq})   is essential. \hfill\Halmos\endproof
}
\end{example}


\proof{Proof of Theorem~\ref{lm:EC}.}
{\rm (i)} We observe that the lower semi-continuity at $s$ of each function $f_n,$ $n\in\N^*,$ follows from lower semi-equicontinuity
of $\{  f_n  \}_{n\in\N^*}$ at $s.$ Thus, to prove statement~(i) it is sufficient to verify (\ref{eq:llim}), which is equivalent to \eqref{eq:P2} because of Remark~\ref{rem:mainconv}.

Let us prove \eqref{eq:P2}. Fix an arbitrary $\varepsilon>0.$
According to Definition~\ref{def:lowerequicont}, there exists $\delta(\varepsilon)>0$
such that for each $n\in\N^*$ and $s^\prime\in B_{\delta(\varepsilon)}(s)$
\begin{equation}\label{eq:P3}
f_n(s^\prime)\ge f_n(s)-\varepsilon.
\end{equation}
Since
\begin{equation}\label{eq:P4}
\ilim_{n\to \infty,s^\prime\rightarrow s}  f_n  (s^\prime)=\sup_{n\ge1,\, \delta>0}\ \  \inf_{k\ge n,\,s^\prime\in B_{\delta}(s)}\ f_k  (s^\prime)\ge \sup_{n\ge1}\ \ \inf_{k\ge n,\,s^\prime\in B_{\delta(\varepsilon)}(s)}\ f_k  (s^\prime),
\end{equation}
\eqref{eq:P3} implies
\begin{equation}\label{eq:P5}
\ilim_{n\to \infty,s^\prime\rightarrow s}  f_n  (s^\prime)\ge \sup_{n\ge1}\inf_{k\ge n}f_k  (s)-\varepsilon=\ilim_{n\to\infty} f_n (s)-\varepsilon,
\end{equation}
where the equalities in (\ref{eq:P4}) and (\ref{eq:P5}) follow from the definition of $\ilim,$
the inequality in (\ref{eq:P4}) holds because
$\{\delta(\varepsilon)\}\subset \{\delta\,:\,\delta>0\},$
and the inequality in (\ref{eq:P5}) follows from (\ref{eq:P3}) and (\ref{eq:P4}). Then,
inequality (\ref{eq:P2}) follows from (\ref{eq:P5}) since $\varepsilon>0$ is arbitrary. Statement~(i) is proved.

{\rm (ii)} We prove statement~(ii) by contradiction. Assume that the sequence of
functions $\{  f_n  \}_{n\in\N^*}$ is not lower semi-equicontinuous at $s.$
Then there exist $\varepsilon^*>0,$ a sequence $\{s_n\}_{n\in\N^*}$ converging to $s,$ and a sequence $\{n_k\}_{k\in\N^*}\subset \N^*$ such that
\begin{equation}\label{eq:eq1}
f_{n_k}(s_k)\le f_{n_k}(s)-\varepsilon^*,\qquad\qquad {   k\in\N^*}.
\end{equation}
  If a sequence $\{n_k\}_{k\in\N^*}$ is bounded (by a positive integer $C$), then (\ref{eq:eq1}) contradicts to lower semi-continuity of each function $f_n,$ $n\in\N^*,C.$ Otherwise, without loss of generality, we may assume that the sequence $\{n_k\}_{k\in\N^*}$ is strictly increasing. Therefore, \eqref{eq:eq1} and \eqref{eq:eq} imply that
\[
\ilim_{n\to \infty,s^\prime\rightarrow s} f_n (s^\prime)\le \lim_{n\to\infty}  f_n (s)-\varepsilon^*.
\]
This is a contradiction to (\ref{eq:llim}). Hence, the sequence of
functions $\{  f_n  \}_{n\in\N^*}$ is lower semi-equicontinuous at $s.$
\hfill\Halmos\endproof

Let us investigate necessary and sufficient conditions for equality (\ref{eq:lim}).

\begin{corollary}\label{cor:limit}
	Let $\{  f_n  \}_{n\in\N^*}$ be a sequence of real-valued functions on a metric
	space $\S$ and $s\in\S.$ If $\{ f_n (s) \}_{n\in\N^*}$ is a convergent sequence, that is, \eqref{eq:eq} holds, then
	the sequence of functions $\{  f_n  \}_{n\in\N^*}$ is equicontinuous at $s$
	if and only if each function $f_n,$ $n\in\N^*,$ is continuous at $s$ and
	equality (\ref{eq:lim}) holds.
\end{corollary}


\proof{Proof.}
Corollary~\ref{cor:limit} follows directly from Theorem~\ref{lm:EC} applied twice to the families $\{ f_n  \}_{n\in\N^*}$ and $\{ - f_n  \}_{n\in\N^*}.$
\hfill\Halmos\endproof


In the following corollary we establish sufficient conditions for lower semi-equicontinuity.


\begin{corollary}[Sufficient conditions for lower semi-equicontinuity]\label{cor:usc-lsec}
Let $\S$ be a metric space and $\{f_n\}_{n\in\N^*}$ be a sequence of
real-valued lower semi-continuous  functions on $\S$  semi-converging uniformly from below to a  real-valued lower semi-continuous function $f$ on $\S.$ If the sequence $\{f_n\}_{n\in\N^*}$ converges pointwise to $f$ on $\S,$ then $\{ f_n \}_{n\in\N^*}$ is lower semi-equicontinuous on $\S.$
\end{corollary}

\proof{Proof.}
If inequality \eqref{eq:P2} holds for all $s\in\S,$ then Remark~\ref{rem:mainconv} and Theorem~\ref{lm:EC}{\rm (ii)} imply that $\{ f_n \}_{n\in\N^*}$ is lower semi-equicontinuous on $\S$ because  the sequence of functions $\{f_n\}_{n\in\N^*}$ converges pointwise to $f$ on $\S.$ Therefore, to finish the proof, let us prove that \eqref{eq:P2} holds  for each $s\in\S.$ Indeed, the uniform semi-convergence from below of $\{f_n\}_{n\in\N^*}$ to $f$ on $\S$ implies that for an arbitrary $\ee > 0$ 
\begin{align}
	\ilim_{n\to\infty, s^\prime\to s} f_n(s^\prime)  \geq f(s) - \ee,
	\label{eq:usc2}
\end{align}
for each $s\in\S.$ Since $\ee>0$ is arbitrarily and $f(s) = \lim_{n\to\infty} f_n(s),$ $s\in\S,$ equality
\eqref{eq:llim} follows from \eqref{eq:P2}. 
\hfill\Halmos\endproof

Let $\S$ be a compact metric space. The Ascoli theorem (see  \cite[p. 96]{HLL96} or \cite[p. 179]{Roy68}) implies that a sequence of
real-valued continuous functions $\{f_n\}_{n\in\N^*}$ on $\S$
converges uniformly on $\S$ to a continuous real-valued function $f$ on $\S$ if and only if $\{f_n\}_{n\in\N^*}$ is equicontinuous and this sequence converges pointwise to $f$ on $\S.$
According to Corollary~\ref{cor:usc-lsec},  a sequence of
real-valued lower semi-continuous functions $\{ f_n \}_{n\in\N^*}$ on $\S,$  converging pointwise to a real-valued lower semi-continuous function $f$ on $\S,$ is lower semi-equicontinuous on $\S$ if $\{f_n\}_{n\in\N^*}$ semi-converges uniformly from below to $f$ on $\S.$ Example~\ref{ex:lsec-usc} illustrates that the converse statement to Corollary~\ref{cor:usc-lsec} does not hold in the general case; that is,
there is a lower semi-equicontinuous sequence $\{ f_n \}_{n\in\N^*}$ of continuous functions on $\S$  converging pointwise to a lower semi-continuous function $f$ such that
$\{ f_n \}_{n\in\N^*}$ does not semi-converge uniformly from below to $f$ on $\S.$


\begin{example}\label{ex:lsec-usc}
{\rm
Define $\S := [0,1]$ endowed with the standard Euclidean metric, $f(s):=\h \{ s\neq 0 \},$ and for $s\in\S$
\begin{align*}
		f_n (s): =
		\begin{cases}
			ns, & \text{if } s\in[0,\frac{1}{n}],  \\
			1, & \text{otherwise. }
		\end{cases}
	\end{align*}
 Then the functions $f_n,$ $n\in\N^*,$ are \textit{continuous} on $\S,$ the function $f$ is \textit{lower semi-continuous} on $\S,$ and the sequence $\{f_n\}_{n\in\N^*}$ \textit{converges pointwise} to $f$ on $\S.$ In addition, the sequence of functions $\{ f_n \}_{n\in\N^*}$
	is \textit{lower semi-equicontinuous} because for each $\ee > 0$ and $s\in\S,$ (i) if $s>0,$ then there exists
	$\delta(s,\ee) = \min\{ s-1/(\lfloor \frac{1}{s} \rfloor + 1), \ee / \lfloor\frac{1}{s}\rfloor \} $
	such that $f_n (s^\prime) \geq f_n (s) - \ee$ for all $n\in\N^*$ and $s^\prime\in B_{\delta(s,\ee)} (s);$
	and (ii) if $s=0,$ then $f_n (s^\prime) \geq 0 = f_n (0)$ for all $n\in\N^*$ and $s^\prime\in\S.$
	The uniform semi-convergence from below of $\{f_n\}_{n\in\N^*}$ to $f$ \textit{does not hold} because $f_n(\frac{1}{n(n+1)})=\frac{1}{n+1}\le 1-\frac12=f(\frac{1}{n(n+1)})-\frac12$ for each $ n\in\N^*,$ that is, the converse statement to Corollary~\ref{cor:usc-lsec} does not hold.~\hfill\Halmos\endproof
}
\end{example}

\section{Fatou's Lemmas in the Classic Form  for Varying Measures}
\label{sec:Fatou}

In this section we establish Fatou's lemmas  in their classic form for varying measures. 
This section   consists of two subsections dealing with weakly and setwise converging measures, respectively.


\subsection{Fatou's lemmas   in  the classic form for weakly converging measures}
\label{subsec:Fatouweak}

%
%
%

The following theorem is the  main result of this subsection .


\begin{theorem}[FL for weakly converging measures]\label{thm:Fatou:w:lsec}
Let $\S$ be a metric space, the sequence of measures $\{\mu_n \}_{n\in\N^*}$ converge  weakly to $\mu\in \M (\S),$
$\{ f_n \}_{n\in\N^*}$ be a lower semi-equicontinuous  sequence of real-valued  functions on $\S,$
and $f$ be a measurable real-valued function on $\S.$
If the following conditions hold:
\begin{itemize}
\item[{\rm (i)}] the sequence $\{f_n\}_{n\in\N^*}$ lower semi-converges to $f$ in measure $\mu;$  
\item[{\rm (ii)}] either $\{ f_n^- \}_{n\in\N^*}$ is a.u.i. w.r.t. $\{ \mu_n \}_{n\in\N^*}$  or Assumption~(ii) of Theorem~\ref{thm:Fatou:w} holds,
\end{itemize}
then
\begin{align}\label{eq:lf:ecMea}
\int_\S f(s) \mu(ds) \leq \ilim\limits_{n\to\infty}\int_\S  f_n (s)\mu_n (ds) .
\end{align}
\end{theorem}

 We recall that asymptotic uniform integrability of $\{ f_n^- \}_{n\in\N^*}$ w.r.t. $\{ \mu_n \}_{n\in\N^*}$ neither implies nor is implied by Assumption~(ii) of Theorem~\ref{thm:Fatou:w}; \cite[Examples~3.1 and 3.2]{FKL18}.
\medskip

\proof{Proof of Theorem~\ref{thm:Fatou:w:lsec}.}
Consider a subsequence $\{ f_{n_k} \}_{k\in\N^*}\subset\{ f_n \}_{n\in\N^*}$ such that
\begin{align}
	 \lim_{k\to\infty} \int_\S f_{n_k} (s) \mu_{n_k} (ds)=\ilim_{n\to\infty}\int_\S  f_n (s)\mu_n (ds).
	\label{eq:subCvg:1}
\end{align}
Assumption {\rm (i)} implies that $\mu ( \{ s\in\S : f_{n_k} (s) \leq f (s) - \ee \} )\to 0$ as $k\to\infty$ for each $\ee > 0.$ Therefore, according to Remark~\ref{rem:pointwise semiconv}, there exists a subsequence $\{ f_{{k_j}} \}_{j\in\N^*}\subset \{ f_{n_k} \}_{k\in\N^*}$ such that $f (s)\le \ilim_{j\to\infty} f_{{k_j}} (s)$ for $\mu$-a.e. $s\in\S. $ Thus, Theorem~\ref{lm:EC}(i) implies that
\[
f (s)\le \ilim_{j\to\infty,s^\prime\to s} f_{{k_j}} (s^\prime),
\]
for $\mu$-a.e. $s\in\S$ and, therefore,
\begin{align}
	\int_\S f (s) \mu (ds) \leq \int_\S \ilim_{j\to\infty,s^\prime\to s} f_{{k_j}} (s^\prime) \mu (ds) .
	\label{eq:subCvg:2}
\end{align}
Theorem~\ref{thm:Fatou:w}, applied to $\{ f_{{k_j}} \}_{j\in\N^*},$ implies
\begin{align}
	\int_\S \ilim_{j\to\infty,s^\prime\to s} f_{{k_j}} (s^\prime)\mu(ds)
	\le \ilim\limits_{j\to\infty}\int_\S  f_{{k_j}} (s)\mu_{{k_j}} (ds) .
	\label{eq:subCvg:3}
\end{align}
Hence, \eqref{eq:lf:ecMea} follows directly from \eqref{eq:subCvg:2}, \eqref{eq:subCvg:3}, and \eqref{eq:subCvg:1}.
\hfill\Halmos\endproof


The following corollary states that the setwise convergence in Theorem~\ref{thm:Fatou:sw} can be substituted by the weak convergence if the integrands form a  lower semi-equicontinuous sequence of functions.

\begin{corollary}[FL for weakly converging measures]\label{cor:Fatou:ec}
Let $\S$ be a metric space, a sequence of measures $\{\mu_n \}_{n\in\N^*}$ converge  weakly to $\mu\in \M (\S),$
$\{ f_n \}_{n\in\N^*}$ be a lower semi-equicontinuous sequence of real-valued functions on $\S.$
If assumption {\rm (ii)} of Theorem~\ref{thm:Fatou:w:lsec} holds,
then  inequality \eqref{eq:ilimn} holds.
\end{corollary}

\proof{Proof.}
Inequality \eqref{eq:ilimn} follows directly from Theorem~\ref{thm:Fatou:w:lsec} and Remark~\ref{rem:pointwise semiconv}.
\hfill\Halmos\endproof


The following example illustrates that Theorem~\ref{thm:Fatou:w:lsec} can provide a more exact
lower bound for the lower limit of integrals than Theorem~\ref{thm:Fatou:w}. 


\begin{example}\label{ex:FatouEC}
{\rm
	Let $\S := [0,2].$ We endow $\S$ with the following metric:
\[
\rho (s_1,s_2) = \h\{ s_1\in[0,1) \}\h\{s_2\in[0,1)\}|s_1-s_2| + \Big(1-\h\{ s_1\in[0,1) \}\h\{s_2\in[0,1)\}\Big)\h\{s_1\neq s_2\}.
\]
To see that $\rho$ is a metric, note that for $s_1,$ $s_2\in\S$ (i) $\rho (s_1,s_2)\in[0,1];$
(ii) $\rho (s_1,s_2)=0$ iff $s_1=s_2;$
(iii) $\rho (s_1,s_2)$ is symmetric in $s_1$ and $s_2;$ and
(iv) for $s_1\neq s_2$ and $s_3\in\S,$ the triangle inequality holds because
$\rho (s_1,s_2) = |s_1-s_2|\leq |s_1-s_3|+|s_3-s_2|= \rho (s_1,s_3)+\rho (s_3,s_2)$ if $s_1,$ $s_2,$ $s_3\in[0,1),$
and $\rho (s_1,s_2)\leq 1\leq \rho (s_1,s_3)+\rho (s_3,s_2)$ otherwise.
	Let $\mu$ be the Lebesgue measure on $\S$ and $\{\mu_n\}_{n\in\N^*}\subset\M(\S)$ be defined as
	\begin{align*}
		\mu_n (C) := \sum_{k=0}^{n-1} \frac{1}{n}I\{\frac{k}{n}\in C\} + \mu (C\cap [1,2]),
		\qquad C\in\SSig, \quad n\in\N^*.
	\end{align*}
	Then the sequence $\{\mu_n\}_{n\in\N^*}$ converges weakly to $\mu$  (see Billingsley~\cite[p. 15, Example 2.2]{Bil68}) and $\{\mu_n\}_{n\in\N^*}$ does not converge setwise to $\mu$  because $\mu_n ([0,1]\setminus\Q) = 0 \not\to 1 = \mu(\Q),$ where $\Q$ is the set of all rational numbers in $[0,1].$
	Define $f \equiv 1$ and
	$f_n (s) = 1-\h \{ s\in (1+\frac{j}{2^k},1+\frac{j+1}{2^k}] \},$ where
	$k = \lfloor \log_2 n \rfloor,$ $j=n-2^k,$ $s\in\S,$ and $n\in\N^*.$
	
	Since the subspace $(1,2]\subset \S$ is endowed with the discrete metric,
	every sequence of functions on $(1,2]$ is equicontinuous.
	Since $f_n (s) = 1$ for $n\in\N^*$ and $s\in[0,1],$ the sequence $\{ f_n \}_{n\in\N^*}$
	is equicontinuous on $[0,1].$ Therefore, $\{ f_n \}_{n\in\N^*}$ is equicontinuous and, thus,
	lower semi-equicontinuous on $\S.$
	In addition, \eqref{eq:sw2aaa:var} holds and $\{ f_n^- \}_{n\in\N^*}$ is a.u.i. w.r.t. $\{\mu_n\}_{n\in\N^*}$
	because $f_n$ is nonnegative for $n\in\N^*.$
	Since $\mu(\{s\in\S:f_n(s) < f(s)\}) = \frac{1}{2^{\lfloor \log_2 n \rfloor}} \to 0$
	as $n\to\infty,$ condition {\rm (i)} from Theorem~\ref{thm:Fatou:w:lsec} holds.
	In view of Theorem~\ref{thm:Fatou:w:lsec},
	\begin{align*}
		& \ilim\limits_{n\to\infty}\int_\S  f_n (s)\mu_n (ds) =
		\ilim\limits_{n\to\infty} \left (\int_0^1  f_n (s)\mu_n (ds) + \int_1^2  f_n (s)\mu_n (ds)\right ) \\
		= & \ilim\limits_{n\to\infty}\left ( 1 + 1  - \frac{1}{2^{\lfloor \log_2 n \rfloor}}\right )  = 2
		\geq \int_\S f(s) \mu(ds) = 2 .
	\end{align*}
	As follows from Theorem~\ref{lm:EC}(i),
	\begin{align*}
	\ilim_{n\to\infty,s^\prime\to s} f_n (s^\prime) = \ilim_{n\to\infty} f_n (s)
	= 1-\h\{s\in[1,2] \} , \qquad s\in\S .
	\end{align*}
	In view of Theorem~\ref{thm:Fatou:w}, \eqref{eq:lf} and \eqref{eq:ilimn} imply
	\begin{align*}
	2 = \ilim\limits_{n\to\infty}\int_\S  f_n (s)\mu_n (ds) \geq \int_\S \ilim_{n\to\infty,s^\prime\to s} f_n (s^\prime) \mu(ds) = \int_\S \ilim_{n\to\infty} f_n (s) \mu(ds) = 1.
	\end{align*}
	Therefore, Theorem~\ref{thm:Fatou:w:lsec} provides a more exact lower bound \eqref{eq:lf:ecMea}
	for the lower limit of integrals than \eqref{eq:lf} and \eqref{eq:ilimn} for weakly
	converging measures and lower semi-equicontinuous sequences of functions.~\hfill\Halmos\endproof
}
\end{example}

\subsection{Fatou's lemmas for setwise converging measures}
\label{subsec:Fatousetwise}

The main results of this subsection, Theorem~\ref{thm:Fatou:sw:ui} and its  Corollary~\ref{thm:Fatou:w:varnotprobsw},  are counterparts for setwise converging measures to Theorem~\ref{thm:Fatou:w:lsec}.


\begin{theorem}[FL for setwise converging measures]\label{thm:Fatou:sw:ui}
Let $(\SS,\SSig)$ be a measurable space, a sequence of measures
$\{\mu_n\}_{n\in\N^*}$ converge setwise to a measure $\mu\in \M(\S),$
and $\{ f_n\}_{n\in\N^*}$ be a sequence of $\olR$-valued measurable functions on $\S.$
If $\{f_n\}_{n\in\N^*}$ lower semi-converges to a real-valued function $f$ in measure $\mu$ and
$\{f_n^-\}_{n\in\N^*}$ is a.u.i. w.r.t. $\{ \mu_n \}_{n\in\N^*},$
then inequality \eqref{eq:lf:ecMea} holds.
\end{theorem}


\proof{Proof.} The proof repeats several lines of the proofs of Theorems~\ref{thm:Fatou:w:lsec} and \ref{thm:Fatou:w}. Consider a subsequence $\{ f_{n_k} \}_{k\in\N^*}\subset\{ f_n \}_{n\in\N^*}$ such that
\begin{align}
	 \lim_{k\to\infty} \int_\S f_{n_k} (s) \mu_{n_k} (ds)=\ilim_{n\to\infty}\int_\S  f_n (s)\mu_n (ds).
	\label{eq:subCvg:1sw}
\end{align}
Since the sequence $\{f_n\}_{n\in\N^*}$ lower semi-converges to $f$ in measure $\mu,$ we have that $\mu ( \{ s\in\S : f_{n_k} (s) \leq f (s) - \ee \} )\to 0$ as $k\to\infty$ for each $\ee > 0.$ Therefore, Remark~\ref{rem:pointwise semiconv} implies that there exists a subsequence $\{ f_{{k_j}} \}_{j\in\N^*}\subset \{ f_{n_k} \}_{k\in\N^*}$ such that $f (s)\le \ilim_{j\to\infty} f_{{k_j}} (s)$ for $\mu$-a.e.~$s\in\S.$ Thus,
\begin{align}
	\int_\S f (s) \mu (ds) \leq \int_\S \ilim_{j\to\infty} f_{{k_j}} (s) \mu (ds) .
	\label{eq:subCvg:2sw}
\end{align}

Now we prove that
\begin{align}
	\int_\S \ilim_{j\to\infty} f_{{k_j}} (s)\mu(ds)
	\le \lim\limits_{j\to\infty}\int_\S  f_{{k_j}} (s)\mu_{{k_j}} (ds) .
	\label{eq:subCvg:3sw}
\end{align}
For this purpose we fix an arbitrary $K > 0.$ Then
\begin{equation}	
\begin{aligned}
	 \ilim_{j\to\infty} \int_\S f_{k_j} (s) \mu_{k_j} (ds) \geq\ilim_{j\to\infty}& \int_\S f_{k_j} (s) \h \{ s\in\S: f_{k_j} (s) > -K \} \mu_{k_j} (ds)\\  +& \ilim_{j\to\infty} \int_\S f_{k_j} (s) \h \{ s\in\S: f_{k_j} (s) \leq -K \} \mu_{k_j} (ds).
	\label{eq:lbd:uisw}
	\end{aligned}
\end{equation}

The following inequality holds:
\begin{equation}\label{eq:finKsw}
\ilim_{j\to\infty} \int_\S f_{k_j} (s) \h \{ s\in\S: f_{k_j} (s) > -K \} \mu_{k_j} (ds)\ge \int_\S \ilim_{j\to\infty}f_{k_j}(s)\mu (ds).
\end{equation}
Indeed, Serfozo's~\cite[Lemma 2.2]{Ser82} applied to the nonnegative sequence $\{f_{k_j} (s) \h \{ s\in\S: f_{k_j} (s) > -K \}+K\}_{j\in\N^*}$ implies
\begin{equation}\label{eq:finK1sw}
\ilim_{j\to\infty} \int_\S f_{k_j} (s) \h \{ s\in\S: f_{k_j} (s) > -K \} \mu_{k_j} (ds)\ge \int_\S \ilim_{j\to\infty}f_{k_j} (s) \h \{ s\in\S: f_{k_j} (s) > -K \}\mu (ds).
\end{equation}
Here we note that
\begin{equation}\label{eq:finK2sw}
f_{k_j}(s)\h \{ s\in\S: f_{k_j}(s)>-K \} \geq f_{k_j}(s),
\end{equation}
for each $s\in\S$ because $K>0.$ Thus, \eqref{eq:finKsw} follows from \eqref{eq:finK1sw} and \eqref{eq:finK2sw}.

Inequalities \eqref{eq:lbd:uisw} and \eqref{eq:finKsw} imply
\[
\begin{aligned}
	 \ilim_{j\to\infty} \int_\S f_{k_j} (s) \mu_{k_j} (ds) \geq \int_\S &\ilim_{j\to\infty,s^\prime\to s}f_{k_j}(s^\prime)\mu (ds)\\  & + \lim_{K\to+\infty}\ilim_{j\to\infty} \int_\S f_{k_j} (s) \h \{ s\in\S: f_{k_j} (s) \leq -K \} \mu_{k_j} (ds),
	\end{aligned}
\]
which is equivalent to \eqref{eq:subCvg:3sw} because $\{ f_{k_j}^- \}_{j\in\N^*}$ is a.u.i. w.r.t. $\{\mu_{k_j}\}_{j\in\N^*}.$

Hence, \eqref{eq:lf:ecMea} follows directly from \eqref{eq:subCvg:2sw}, \eqref{eq:subCvg:3sw}, and \eqref{eq:subCvg:1sw}.
\hfill\Halmos\endproof


The following corollary to Theorem~\ref{thm:Fatou:sw:ui}  generalizes
Theorem~\ref{thm:Fatou:sw}.


\begin{corollary}\label{thm:Fatou:w:varnotprobsw}
Let $(\SS,\SSig)$ be a measurable space, a sequence of measures
$\{\mu_n\}_{n\in\N^*}$ converge setwise to a measure $\mu\in \M(\S),$
and $\{ f_n\}_{n\in\N^*}$ be a sequence of $\olR$-valued measurable functions on $\S$
 lower semi-converging to a real-valued function $f$ in measure $\mu.$ If there exists a sequence of measurable real-valued functions $\{g_n\}_{n\in\N^*}$
	on $\SS$ such that $ f_n (s)\ge g_n(s)$ for all $n\in\N^*$ and $s\in \SS,$ and if
\eqref{eq:sw2aaa:sw:var} holds,
then inequality \eqref{eq:lf:ecMea} holds.	
\end{corollary}


\proof{Proof.} Consider an increasing sequence $\{ n_k \}_{k\in\N^*}$ of natural numbers such that
\begin{align}
	 \lim_{k\to \infty} \int_\S f_{n_k} (s) \mu_{n_k} (ds)=\ilim_{n\to\infty}\int_\S  f_n (s)\mu_n (ds).
	\label{eq:subCvg:1sw1}
\end{align}
Since the sequence $\{f_n\}_{n\in\N^*}$ lower semi-converges to $f$ in measure $\mu,$ we have that $\mu ( \{ s\in\S : f_{n_k} (s) \leq f (s) - \ee \} )\to 0$ as $k\to\infty$ for each $\ee > 0.$ Therefore, Remark~\ref{rem:pointwise semiconv} implies that there exists a subsequence $\{ f_{{k_j}} \}_{j\in\N^*}\subset \{ f_{n_k} \}_{k\in\N^*}$ such that $f (s)\le \ilim_{j\to\infty} f_{{k_j}} (s)$ for $\mu$-a.e.~$s\in\S. $ Thus,
\begin{align}
	\int_\S f (s) \mu (ds) \leq \int_\S \ilim_{j\to\infty} f_{{k_j}} (s) \mu (ds) .
	\label{eq:subCvg:2sw2}
\end{align}

Now we prove that
\begin{align}
	\int_\S \ilim_{j\to\infty} f_{{k_j}} (s)\mu(ds)
	\le \lim\limits_{j\to\infty}\int_\S  f_{{k_j}} (s)\mu_{{k_j}} (ds) .
	\label{eq:subCvg:3sw3}
\end{align}
Theorem~\ref{thm:Fatou:sw:ui}, applied to the sequence $\{ f_{k_j}-g_{k_j} \}_{j\in\N^*},$ implies
\[
\begin{aligned}
&\int_\S \ilim\limits_{j\to\infty} f_{k_j} (s)\mu(ds)-\int_\S \slim\limits_{j\to\infty} g_{k_j} (s)\mu(ds) \\
& \le \int_\S \ilim\limits_{j\to\infty} (f_{k_j} (s)-g_{k_j}(s))\mu(ds)\le \ilim\limits_{j\to\infty}\int_\S f_{k_j} (s)\mu_{k_j} (ds)-\slim\limits_{j\to\infty}\int_\S g_{k_j} (s)\mu_{k_j} (ds),
	\end{aligned}
\]
where the first and third inequalities follow from the basic properties of infimums and supremums. Then, \eqref{eq:sw2aaa:sw:var} implies \eqref{eq:subCvg:3sw3}.
Hence, \eqref{eq:lf:ecMea} follows directly from \eqref{eq:subCvg:2sw2}, \eqref{eq:subCvg:3sw3}, and \eqref{eq:subCvg:1sw1}.
\hfill\Halmos\endproof


Theorem~\ref{thm:Fatou:sw:ui} provides a more exact lower bound for the lower limit of integrals than
Theorem~\ref{thm:Fatou:sw}. This fact is illustrated in Example~\ref{ex:Fatou:sw}.
	

\begin{example}[{cp. Feinberg et al.~\cite[Example 4.1]{FKZ16}}]\label{ex:Fatou:sw}
{\rm
	Let $\SS = [0,1],$ $\SSig=\B([0,1]),$ $\mu$ be the Lebesgue measure on $\S,$ and for $C\in\B(\SS)$
	\begin{align*}
		\mu_n (C) := \int_C 2\h\{s\in\SS :\frac{2k}{2^n}<s<\frac{2k+1}{2^n},k=0,1,2,\ldots,2^{n-1}-1\}\mu(ds), \quad n\in\N^*.
	\end{align*}
	Define $f \equiv 1$ and
	$f_n (s) = 1-\h \{ s\in [\frac{j}{2^k},\frac{j+1}{2^k}] \},$ where
	$k = \lfloor \log_2 n \rfloor,$ $j=n-2^k,$ $s\in\S,$ and $n\in\N^*.$
	Then the sequence $\{\mu_n\}_{n\in\N^*}$ converges setwise to $\mu,$ \eqref{eq:sw2aaa:sw:var} holds and $\{ f_n^- \}_{n\in\N^*}$ is a.u.i. w.r.t. $\{ \mu_n \}_{n\in\N^*},$ and
	the sequence $\{f_n\}_{n\in\N^*}$ lower semi-converges to $f$ in measure $\mu.$
	In view of Theorem~\ref{thm:Fatou:sw:ui} and \eqref{eq:ilimn},
	\begin{align*}
		1 = \ilim\limits_{n\to\infty}\int_\SS  f_n (s)\mu_n (ds) \geq \int_\SS f(s) \mu(ds) = 1>0=\int_\SS \ilim_{n\to\infty} f_n (s) \mu(ds) .
	\end{align*}
	Therefore, Theorem~\ref{thm:Fatou:sw:ui} provides a more exact lower bound
	for the lower limit of integrals than inequality \eqref{eq:ilimn}. 
	$\ $ \hfill\Halmos\endproof
}
\end{example}	
%
%
%
%
%

\section{Lebesgue's Convergence Theorem for Varying Measures}
\label{sec:main}

In this section, we present Lebesgue's convergence
theorem for varying measures $\{\mu_n\}_{n\in\N^*}$ and functions that
are a.u.i.~w.r.t. $\{\mu_n\}_{n\in\N^*}.$
The following corollary follows from Theorem~\ref{thm:Fatou:w}.
It also follows from Serfozo~\cite[Theorem 3.5]{Ser82} adapted to general metric spaces.
We provide it here for completeness.


\begin{corollary}[{Lebesgue's convergence theorem for weakly converging measures~\cite[Corollary~2.8]{FKL18}}]\label{cor:leb:weak}
Let $\S$ be a metric space, $\{\mu_n\}_{n\in\N^*}$ be a sequence of measures on $\S$  converging weakly to
$\mu\in \M (\S),$ and $\{ f_n \}_{n\in\N^*}$ be an a.u.i. (see \eqref{eq:tv:ui:lim})
w.r.t. $\{\mu_n\}_{n\in\N^*}$ sequence of measurable $\olR$-valued functions on $\S$ such that
$\lim_{n\to \infty,\, s'\to s} f_n (s')$ exists for $\mu$-a.e. $s\in\S,$
then
\begin{align*}
	\lim\limits_{n\to\infty}\int_\S  f_n (s)\mu_n (ds) = \int_\S \lim_{n\to\infty,\, s'\to s} f_n (s')\mu(ds) = \int_\S \lim_{n\to\infty} f_n (s)\mu(ds) .
\end{align*}
\end{corollary}




The following corollary states the  convergence theorem for weakly converging
measures $\mu_n$ and for an equicontinuous sequence of functions $\{ f_n \}_{n\in\N^*}$.


\begin{corollary}[Lebesgue's convergence theorem for weakly converging measures]\label{cor:leb:weak:eq}
Let $\S$ be a metric space, the sequence of measures $\{\mu_n \}_{n\in\N^*}$ converge weakly to $\mu\in \M (\S),$
$\{ f_n \}_{n\in\N^*}$ be a sequence of real-valued equicontinuous functions on $\S,$
and $f$ be a measurable real-valued function on $\S.$
If the sequence $\{f_n\}_{n\in\N^*}$ converges to $f$ in measure $\mu$ and is
a.u.i. (see \eqref{eq:tv:ui:lim}) w.r.t. $\{\mu_n\}_{n\in\N^*},$ then
\begin{align}
	\lim\limits_{n\to\infty}\int_\S  f_n (s)\mu_n (ds) = \int_\S f(s) \mu(ds) .
	\label{eq:leb:cvg}
\end{align}
\end{corollary}


\proof{Proof.}
Corollary~\ref{cor:leb:weak:eq} follows from Theorem~\ref{thm:Fatou:w:lsec} applied to $\{ f_n \}_{n\in\N^*}$ and  $\{- f_n \}_{n\in\N^*}.$ \hfill\Halmos\endproof

The following corollary for setwise converging measures follows directly from Theorem~\ref{thm:Fatou:sw:ui}.
	

\begin{corollary}[Lebesgue's convergence theorem for setwise converging measures]\label{cor:dc:sw:ui}
Let $(\SS,\SSig)$ be a measurable space, a sequence of measures
$\{\mu_n\}_{n\in\N^*}$ converge setwise to a measure $\mu\in \M(\S),$
and $\{ f_n \}_{n\in\N^*}$ be a sequence of $\olR$-valued measurable functions on $\S.$
If the sequence $\{f_n\}_{n\in\N^*}$ converges to a measurable real-valued function $f$ in measure $\mu$ and this sequence is
a.u.i. (see \eqref{eq:tv:ui:lim}) w.r.t. $\{\mu_n\}_{n\in\N^*},$
then \eqref{eq:leb:cvg} holds.
\end{corollary}


\proof{Proof.}
Corollary~\ref{cor:dc:sw:ui} follows from Theorem~\ref{thm:Fatou:sw:ui}  applied to 
$\{ f_n \}_{n\in\N^*}$ and  $\{- f_n \}_{n\in\N^*}.$
\hfill\Halmos\endproof
%
%

\section{Monotone Convergence Theorem for Varying Measures}
\label{sec:swtv}

In this section, we present monotone convergence
theorems for varying measures. 


\begin{theorem}[Monotone convergence theorem for weakly converging measures]\label{thm:mono}
Let $\S$ be a metric space, $\{\mu_n\}_{n\in\N^*}$ be a sequence of measures on $\S$ that converges weakly to
$\mu\in \M (\S),$ $\{ f_n \}_{n\in\N^*}$ be a sequence of lower semi-continuous $\olR$-valued functions on $\S$ such that $f_n(s)\le f_{n+1}(s)$ for each $n\in\N^*$ and $s\in\S,$ and $f(s) := \lim_{n\to \infty} f_n (s),$ $s\in\S.$ If the following conditions hold:
	\begin{itemize}
		\item[{\rm(i)}]the function $f$ is upper semi-continuous;
		\item[{\rm(ii)}] the functions $f_1^-$ and $f^+$ are a.u.i. w.r.t. $\{\mu_n\}_{n\in\N^*};$
	\end{itemize}
	then \eqref{eq:leb:cvg} holds.
\end{theorem}

\begin{remark}\label{rm:mono}
The lower semi-continuity of $f_n$ and pointwise convergence of $f_n$ to $f$ imply
the lower semi-continuity of $f.$ Therefore, under assumptions in Theorem~\ref{thm:mono} the function $f$
is continuous.
\end{remark}


The following example demonstrates the necessity of the condition {\rm (i)} in
Theorem~\ref{thm:mono}.


\begin{example}\label{ex:mono:ii}
{\rm
	Consider $\S = [0,1]$ endowed with the standard Euclidean metric,
	$f (s) = \h \{s\in(0,1]\},$ $s\in\S,$ $f_n (s) = \min\{ n s, 1 \},$
	$n\in\N^*$ and $s\in\S,$ and probability measures
	\begin{align}
		\mu_n (C) := \int_C n\h \{ s\in[0,\frac{1}{n}] \}\nu(ds), \quad
		\mu (C) := \h \{ 0 \in C \}, \qquad {  C\in\B(\SS),} \quad n\in\N^*,
		\label{eq:ex:pm}
	\end{align}
where $\nu$ is the Lebesgue measure on $\SS.$

Then $f_n (s) \uparrow f (s)$ for each $s\in\S$
	as $n\to\infty,$ and the sequence of probability measures $\mu_n$ converges weakly to $\mu.$
	Since the functions $f_1$ and $f$ are bounded, condition {\rm (ii)} from
	Theorem~\ref{thm:mono} holds.
	The function $f_n$ is continuous, and the function $f$ is lower semi-continuous,
	but $f$ is not upper semi-continuous.
	Since $\int_\S f_n (s)\mu_n (ds) = \frac{1}{2},$ $n\in\N^*,$ and $\int_\S f(s) \mu(ds) = 0,$
	formula \eqref{eq:leb:cvg} does not hold.
	\hfill\Halmos\endproof
}
\end{example}


\proof{Proof of Theorem~\ref{thm:mono}.}
Since $f_n (s) \leq f(s),$
\begin{align*}
	f (s)=\ilim_{n\to\infty,s'\to s} f_n (s') \leq \slim_{n\to\infty,s'\to s} f_n (s') \leq
	\slim_{s'\to s} f(s') \leq f (s),\qquad s\in\S,
\end{align*}
where the first equality follows from Lemma~\ref{lm:mono:lim}, and the last inequality holds because $f$ is upper semi-continuous. Hence, $\lim_{n\to\infty,s'\to s} f_n (s') = f (s),$ $s\in\S.$
In addition, condition {\rm (ii)} implies that the sequence $\{f_n\}_n$ is
a.u.i. w.r.t. $\{\mu_n\}_{n\in\N^*}.$ Therefore, Corollary~\ref{cor:leb:weak} implies \eqref{eq:leb:cvg}.
\hfill\Halmos\endproof


\begin{corollary}\label{cor:mono}
	Let $\S$ be a metric space, $\{\mu_n\}_{n\in\N^*}$ be a sequence of measures on $\S$ that converges weakly to
$\mu\in \M (\S),$ and $\{ f_n \}_{n\in\N^*}$ be a pointwise nondecreasing sequence of measurable $\olR$-valued functions on $\S.$ Let $f (s) := \lim_{n\to \infty} f_n (s)$ and $\underline{f}_n (s) := \ilim_{s^\prime\to s} f_n (s^\prime) ,$ $s\in\S.$ If
	\begin{itemize}
\item[{\rm(i)}] the function $f$ is real-valued and upper semi-continuous;
		\item[{\rm(ii)}] the sequence $\{\underline{f}_n\}_{n\in\N^*}$ lower semi-converges to $f$ in measure $\mu;$ and
		\item[{\rm(iii)}] the functions $\underline{f}_1^-$ and $f^+$ are a.u.i. w.r.t. $\{\mu_n\}_{n\in\N^*};$
	\end{itemize}
	then \eqref{eq:leb:cvg} holds.
\end{corollary}


The following example demonstrates the necessity of condition {\rm (ii)} from
Corollary~\ref{cor:mono}.


\begin{example}\label{ex:mono:i}
{\rm
	Consider $\S = [0,1]$ endowed with the standard Euclidean metric, $f (s) = 1,$
	\begin{align*}
		f_n (s) =
		\begin{cases}
			1, & \text{if } s = 0,  \\
			\min\{ n s, 1 \}, & \text{if } s\in (0,1],
		\end{cases}
		\qquad n\in\N^*, \quad s\in\S,
	\end{align*}
	and probability measures $\mu_n,$ $n\in\N^*,$ and $\mu$ defined in \eqref{eq:ex:pm}.
	Then $\underline{f}_n (s) = \min\{ n s, 1 \}, $  $f_n (s) \uparrow f (s)$
	for each $s\in\S$ as $n\to\infty,$ and the sequence of probability measures $\mu_n$
	converges weakly to $\mu.$ Since the functions $\underline{f}_1$ and $f$ are bounded, condition {\rm (iii)} from
	Corollary~\ref{cor:mono} holds.
	Condition {\rm (ii)} from Corollary~\ref{cor:mono} does not hold because $f(0)=f_n(0) =1$
	and $\underline{f}_n (0) = 0$ for each $n\in\N^*.$
	Since $\int_\S f_n (s)\mu_n (ds) = \frac{1}{2},$ $n\in\N^*,$
	and $\int_\S f(s) \mu(ds) = 1,$ formula \eqref{eq:leb:cvg} does not hold.
	\hfill\Halmos\endproof
}
\end{example}


\proof{Proof of Corollary~\ref{cor:mono}.}
Since the function $\underline{f}_n$ is lower semi-continuous, Theorem~\ref{thm:mono} implies
\begin{align}
	\lim\limits_{n\to\infty}\int_\S  \underline{f}_n (s)\mu_n (ds) = \int_\S \lim\limits_{n\to\infty} \underline{f}_n(s) \mu(ds) .
	\label{eq:leb:cvg:lower}
\end{align}
Condition {\rm (i)} implies that there exists a subsequence $\{ f_{n_k} \}_{k\in\N^*}\subset\{ f_n \}_{n\in\N^*}$ such that
\begin{align}
	\ilim\limits_{k\to\infty} \underline{f}_{n_k} (s) \geq f (s) \quad \text{for }\mu\text{-a.e. } s\in \S .
	\label{eq:m:cvae}
\end{align}	
Since $\underline{f}_n (s)\leq f_n(s)\leq f(s),$ $n\in\N^*$ and $s\in\S,$ and the sequence
$\{ \underline{f}_n \}_{n\in\N^*}$ is pointwise nondecreasing, \eqref{eq:m:cvae} implies that
\begin{align}
	f(s) = \lim_{n\to\infty} \underline{f}_n (s) \quad \text{for } \mu\text{-a.e. } s\in\S .
	\label{eq:m:cvae1}
\end{align}
Hence, \eqref{eq:leb:cvg:lower} and \eqref{eq:m:cvae1} imply
\begin{align}
	\lim\limits_{n\to\infty}\int_\S  \underline{f}_n (s)\mu_n (ds) = \int_\S f (s) \mu(ds) .
	\label{eq:leb:cvg:lower:1}
\end{align}
Since $\underline{f}_n (s)\leq f_n(s)\leq f(s),$ $n\in\N^*$ and $s\in\S,$
\begin{align}
	\lim\limits_{n\to\infty}\int_\S  \underline{f}_n (s)\mu_n (ds) \leq \ilim\limits_{n\to\infty}\int_\S  f_n (s)\mu_n (ds) \leq \slim\limits_{n\to\infty}\int_\S  f_n (s)\mu_n (ds)  \leq \slim\limits_{n\to\infty} \int_\S f (s) \mu_n(ds) .
	\label{eq:leb:cvg:lower:2}
\end{align}
Theorem~\ref{thm:Fatou:w} applied to the sequence $\{ -f\}$ and the upper semi-continuity of $f$ imply
\begin{align}
	\slim\limits_{n\to\infty} \int_\S f (s) \mu_n(ds) \leq \int_\S f (s) \mu(ds).
	\label{eq:leb:cvg:lower:3}
\end{align}
Therefore, \eqref{eq:leb:cvg:lower:1}, \eqref{eq:leb:cvg:lower:2}, and \eqref{eq:leb:cvg:lower:3} imply \eqref{eq:leb:cvg}.
\hfill\Halmos\endproof

The following corollary from Theorem~\ref{thm:Fatou:sw:ui} is the counterpart to Theorem~\ref{thm:mono}
for setwise converging measures.


\begin{corollary}[Monotone convergence theorem for setwise converging measures]\label{cor:mono:sw}
Let $(\SS,\SSig)$ be a measurable space, a sequence of measures
$\{\mu_n\}_{n\in\N^*}$ converge setwise to a measure $\mu\in \M(\S),$
and $\{ f_n \}_{n\in\N^*}$ be a pointwise nondecreasing sequence of
	measurable $\olR$-valued functions on $\SS.$ Let $f (s) := \lim_{n\to \infty} f_n (s),$ $s\in\SS.$
	If the functions $f_1^-$ and $f^+$ are a.u.i. w.r.t. $\{\mu_n\}_{n\in\N^*},$	then \eqref{eq:leb:cvg} holds.
\end{corollary}


\proof{Proof.}
Since $f_n\uparrow f,$ \eqref{eq:leb:cvg} follows directly from Theorem~\ref{thm:Fatou:sw:ui}  applied to the sequences $\{ f_n\}_{n\in\N^*}$ and $\{ -f_n\}_{n\in\N^*}.$
\hfill\Halmos\endproof

\section{Applications to Markov Decision Processes}
\label{sec:appl}

Consider a discrete-time MDP with a state space $\X,$ an action space $\A,$ one-step
costs $c,$ and transition probabilities $q.$ Assume that $\X$ and $\A$ are Borel subsets
of Polish (complete separable metric) spaces. 
Let  $c(x,a):\X\times\A\mapsto\olR$ be the one-step cost and $q(B|x,a)$ be the transition kernel
representing the probability that the next state is in $B\in\mathcal{B}(\X),$
given that the action $a$ is chosen at the state $x.$  The cost function $c$ is assumed to be measurable and bounded below.

The decision process proceeds as follows: at each time epoch
$t=0,1,\dots,$ the current state of the system, $x,$ is observed.
A decision-maker chooses an action $a,$ the cost $c(x,a)$ is
accrued, and the system moves to the next state according to
$q(\cdot|x,a).$ Let $H_t = (\X\times\A)^{t}\times\X$
be the set of histories for $t=0,1,\dots\ .$ A
(randomized) decision rule at period $t=0,1,\dots$ is a regular transition probability
$\pi_t : H_t\mapsto \A,$ 
that is, (i) $\pi_t(\cdot|h_t)$ is a probability distribution on $\A,$
where $h_t=(x_0,a_0,x_1,\dots,a_{t-1},x_t),$ and (ii) for any measurable subset
$B\subset \A,$ the function $\pi_t(B|\cdot)$ is measurable on $H_t.$
A policy $\pi$ is a sequence $(\pi_0,\pi_1,\dots)$ of decision rules.  Let $\PS$ be the set of all policies.
A policy $\pi$ is called non-randomized if each probability measure $\pi_t(\cdot|h_t)$ is
concentrated at one point. A non-randomized policy is called stationary if all decisions depend
only on the current state.

The Ionescu Tulcea theorem implies that an initial state $x$ and a policy $\pi$ define a unique
probability $\PR_{x}^{\pi}$ on the set of all trajectories $\mathbb{H}_{\infty}=(\X\times\A)^{\infty}$ endowed with the product of $\sigma$-fields defined by Borel $\sigma$-fields of $\X$ and $\A;$ see Bertsekas and Shreve~\cite[pp. 140--141]{BS96} or Hern\'{a}ndez-Lerma and Lasserre~\cite[p. 178]{HLL96}.
Let $\mathbb{E}_{x}^{\pi}$ be an expectation w.r.t. $\PR_{x}^{\pi}.$

For a finite-horizon $N\in\N^*,$ let us define the expected total discounted costs,
\begin{align}\label{eqn:sec_model def:finite total disc cost}
    v_{N,\a}^{\pi} (x):= \mathbb{E}_{x}^{\pi} \sum_{t=0}^{N-1}
    \alpha^{t} c(x_t,a_t), \;\; x\in\X,
\end{align}
where $\alpha\in [0,1]$ is the discount factor.
When $N=\infty$ and $\alpha\in [0,1),$
equation (\ref{eqn:sec_model def:finite total disc cost}) defines an
infinite-horizon expected total discounted cost denoted by $v_{\a}^{\pi}(x).$ Let $v_\alpha (x):=\inf_{\pi\in\PS} v_\a^\pi(x),$ $x\in\X.$  A policy $\pi$ is called optimal for the discount factor $\a$ if $v^\pi_\a(x)=v_\a(x)$ for all $x\in\X.$

The \emph{average cost per unit time} is defined as
\begin{align*}
    w_1^{\pi}(x):=\limsup_{N\to \infty} \frac{1}{N}v_{N,1}^{\pi} (x), \;\; x\in\X.
\end{align*}
Define the optimal value function $w_1(x):=\inf_{\pi\in\Pi} w_1^{\pi} (x),$  $x\in\X.$
A policy $\pi$ is called average-cost optimal if $w_1^{\pi}(x)=w_1(x)$ for all $x\in\X.$

We remark that in general action sets may depend on current states, and usually the state-dependent sets $A(x)$ are considered for all $x\in\X$.  In our problem formulations $A(x)=\A$ for all $x\in\X.$  This problem formulation is simpler than a  formulation with the sets $A(x),$ and these two problem formulations are equivalent because we allow that $c(x,a)=+\infty$ for some $(x,a)\in\X\times\A.$  For example, we may set $A(x)=\{a\in\A: c(x,a)<+\infty\}.$  For  a  formulation with the sets $A(x),$ one may define $c(x,a)=+\infty$ when $a\in\A \setminus A(x)$ and use the action sets $\A$ instead of $A(x).$

To establish the existence of the average-cost optimal policies via an optimality inequality for problems with compact action sets,
Sch\"{a}l~\cite{Sch93} considered two continuity conditions \textbf{W} and \textbf{S} for problems with weakly and setwise continuous transition probabilities, respectively.
For setwise continuous transition probabilities,
Hern\'{a}ndez-Lerma~\cite{HL91} generalized
Assumption~\textbf{S} to Assumption~\textbf{S*}  to cover MDPs with possibly noncompact action sets.
For the similar purpose, when transition probabilities are weakly continuous,
Feinberg~et~al.~\cite{FKZ12} generalized Assumption~\textbf{W} to Assumption~\textbf{W*}.

We recall that a function $f:\U\mapsto \olR$ defined on a metric space $\U$  is
called inf-compact (on $\U$), if for every $\lambda\in\R$ the level set $\{u\in\U :f(u)\leq \lambda \}$
is compact. A subset of a metric space is also a metric space with respect to the same metric.  For $U\subset \U,$ if the domain of $f$ is narrowed to $U,$ then this function is  called the restriction of $f$ to $U.$
\begin{definition}[{Feinberg et al.~\cite[Definition 1.1]{FKZ13}, Feinberg~\cite[Definition 2.1]{Ftut}}]
\label{def:k inf compact}
	A function $f:\X\times\A\mapsto \olR$ is called $\K$-inf-compact, if for
	every nonempty compact subset $\mathcal{K}$ of $\X$ the restriction of $f$ to $\mathcal{K}\times\A$ is an inf-compact function.
\end{definition}

\noindent
\textbf{Assumption W*} ({\rm Feinberg~et~al.~\cite{FKZ12, POMDP}}, Feinberg and Lewis~\cite{FLe17}, or Feinberg~\cite{Ftut})\textbf{.}

(i) the function $c$ is $\K$-inf-compact; 

(ii) the transition probability $q(\cdot|x,a)$ is weakly continuous in
  $(x,a)\in \X\times\A.$

\vspace{0.1in}

\noindent
\textbf{Assumption S*} ({\rm Hern\'{a}ndez-Lerma~\cite[Assumption 2.1]{HL91}})

(i) the function $c(x,a)$ is inf-compact in $a\in \A$ for each $x\in\X;$

(ii) the transition probability $q(\cdot|x,a)$ is setwise continuous in
  $a\in \A$ for each $x\in\X .$
\vspace{0.1in}


Let
\begin{align}\label{defmauaw}
	\begin{split}
  		m_{\alpha}: = \underset{x\in\X}{\inf} v_{\alpha}(x), & \quad
  		u_{\alpha}(x): = v_{\alpha}(x) - m_{\alpha}, \\
  		\underline{w}: = \underset{\alpha\uparrow1}{\liminf}(1-\alpha)m_{\alpha}, & \quad
  		\bar{w}: = \underset{\alpha\uparrow1}{\limsup}(1-\alpha)m_{\alpha} .
	\end{split}
\end{align}
The function $u_\alpha$ is called the discounted relative value function.  If either Assumption \textbf{W*} or Assumption \textbf{S*} holds, let us consider the following assumption.

\vspace{0.1in}
\noindent
\textbf{Assumption B.}
(i) $w^{*} := \inf_{x\in\X} w_1 (x)< +\infty;$ and
(ii) $\underset{\alpha\in [0,1)}{\sup} u_{\alpha}(x) < +\infty,$ $x\in\X.$

\vspace{0.1in}

As follows from Sch\"al~\cite[Lemma 1.2(a)]{Sch93},
Assumption~\textbf{B}(i) implies that $m_\alpha< +\infty$
for all $\alpha\in [0,1).$  Thus, all the quantities in
\eqref{defmauaw} are defined.  

It is known \cite{FKZ12,Sch93} that, if a stationary policy $\phi$ satisfies the average-cost optimality inequality (ACOI)
\begin{align}
	\underline{w} + u(x) &\geq c(x,\phi (x)) + \int_{\X} u(y)q(dy|x,\phi (x)) ,
    \qquad x\in\X,
	\label{eqn:ACOI}
\end{align}
for some nonnegative measurable function $u:\X\to{\bf R},$ then the stationary policy $\phi$ is average-cost optimal.
A nonnegative measurable function $u(x)$ satisfying inequality \eqref{eqn:ACOI} with some stationary policy $\phi$ is
called an average-cost relative value function. The following two theorems state the validity of
the ACOI  under Assumptions~\textbf{W*}
(or Assumption~\textbf{S*}) and \textbf{B}.

\begin{theorem}[{Feinberg et al.~\cite[Corollary 2 and p. 603]{FKZ12}}]\label{thm:ACOI:w}
Let Assumptions \textbf{W*} and \textbf{B} hold. For an arbitrary sequence $\{ \a_n\uparrow 1 \}_{n\in\N^*},$
let
\begin{align}
  	u(x) := \liminf_{n\to \infty,y\rightarrow x} u_{\a_n}(y) , \qquad x\in\X .
  	\label{EQN1}
\end{align}
Then there exists a  stationary  policy $\phi$ satisfying ACOI \eqref{eqn:ACOI} with the function $u$ defined in
\eqref{EQN1}.  Therefore, $\phi$ is a stationary average-cost optimal policy.  In addition, the function $u$ is lower semi-continuous, and

\begin{align}\label{eqsfkz}
    w_1^{\phi}(x)= \underline{w} = \lim_{\alpha\uparrow 1}
    (1-\alpha)v_{\alpha} (x)= \lim_{\alpha\uparrow 1}
    (1-\alpha)m_{\alpha}  = \bar{w} = w^*, \qquad x\in\X.
\end{align}
\end{theorem}

\begin{theorem}[{Hern\'{a}ndez-Lerma~\cite[Section 4, Theorem]{HL91}}]\label{thm:ACOI:s}
Let Assumptions \textbf{S*} and \textbf{B} hold.  For an arbitrary sequence $\{ \a_n\uparrow 1 \}_{n\in\N^*},$
let
\begin{align}
  	u(x) := \liminf_{n\to \infty} u_{\a_n}(x) , \qquad x\in\X .
  	\label{eqn:tu:setwise}
\end{align}
Then there exists a  stationary  policy $\phi$ satisfying ACOI \eqref{eqn:ACOI} with the function $u$ defined in
\eqref{eqn:tu:setwise}.  Therefore, $\phi$ is a stationary average-cost optimal policy.  In addition,   \eqref{eqsfkz} holds.
\end{theorem}

The following corollary from Theorem~\ref{thm:ACOI:w} provides a sufficient condition for the validity of ACOI~\eqref{eqn:ACOI} with a relative value function $u$ defined in \eqref{eqn:tu:setwise}.
\begin{corollary}\label{c6.4EF} Let Assumptions \textbf{W*} and \textbf{B} hold and there exist  a sequence $\{ \a_n\uparrow1 \}_{n\in\N^*}$ of nonnegative discount factors such that
 the sequence of functions $\{u_{\a_n}\}_{n\in\N^*}$ is lower semi-equicontinuous.  Then the conclusions of Theorem~\ref{thm:ACOI:w} hold for the function $u$ defined in \eqref{eqn:tu:setwise} for this sequence $\{ \a_n\}_{n\in\N^*}.$
\end{corollary}
\proof{Proof.} Since the sequence of functions $\{u_{\a_n}\}_{n\in\N^*}$ is lower semi-equicontinuous, the functions $u$ defined in \eqref{EQN1} and in \eqref{eqn:tu:setwise} coincide in view of Theorem~\ref{lm:EC}(i).
\hfill\Halmos\endproof

Consider the following equicontinuity condition (EC) on the discounted relative value functions.

\vspace{0.1in}
\noindent
\textbf{Assumption EC.} There exists a sequence $\{ \a_n \}_{n\in\N^*}$ of nonnegative discount factors such that $\a_n\uparrow1$ as $n\to\infty,$ and the following two conditions hold:

{\rm (i)} the sequence of functions $\{u_{\a_n}\}_{n\in\N^*}$ is equicontinuous;

{\rm (ii)} there exists a nonnegative measurable function $U(x),$ $x\in\X,$ such that
$U(x)\geq u_{\a_n}(x),$ $n\in\N^*,$ and $\int_{\X} U(y)q(dy|x,a) < +\infty$ for all $x\in\X$ and $a\in \A.$
\vspace{0.1in}

Under each of the Assumptions \textbf{W*} or \cite[Assumption 4.2.1]{HLL96}, which is stronger than \textbf{S*}, and under Assumptions \textbf{B} and \textbf{EC}, there exist a sequence $\{ \a_n\uparrow 1 \}_{n\in\N^*}$ of nonnegative discount factors and a stationary policy $\phi$ satisfying the average-cost optimality equations (ACOEs)
\begin{align}
	w^* + u(x) = c(x,\phi(x)) + \int_{\X} u(y)q(dy|x,\phi(x)) =
  		\min_{a\in \A} \left[ c(x,a) + \int_{\X} u(y)q(dy|x,a) \right],
\label{eqn:ACOE}
\end{align}
with $u$ defined in \eqref{EQN1} for the sequence $\{ \a_n\uparrow 1 \}_{n\in\N^*},$ and the function $u$ is continuous; see Feinberg and Liang~\cite[Theorem~3.2]{FLi17} for \textbf{W*} and
Hern\'{a}ndez-Lerma and Lasserre~\cite[Theorem 5.5.4]{HLL96}.
We remark that the quantity $w^*$ in \eqref{eqn:ACOE} can be replaced with any other quantity in \eqref{eqsfkz}.

In addition, since the
left equation in \eqref{eqn:ACOE} implies inequality \eqref{eqn:ACOI}, every stationary policy $\phi$ satisfying
\eqref{eqn:ACOE} is average-cost optimal.  Observe that in these cases the function $u$ is continuous (see \cite[Theorem~3.2]{FLi17} for \textbf{W*} and \cite[Theorem 5.5.4]{HLL96}), while under conditions of Theorems~\ref{thm:ACOI:w} and \ref{thm:ACOI:s} the corresponding functions $u$ may not be continuous; see Examples~\ref{ex:mdpLEC}, \ref{ex:mdpLEC:s}.
Below we provide more general conditions for the validity of the ACOEs.  In particular, under these conditions the relative value functions $u$ may not be continuous.

Now, we introduce Assumption~\textbf{LEC}, which is weaker than Assumption~\textbf{EC}.  Indeed, Assumption~\textbf{EC}(i) is obviously stronger than \textbf{LEC}(i). In view of the Ascoli  theorem (see  \cite[p. 96]{HLL96} or \cite[p. 179]{Roy68}), \textbf{EC}(i) and the first claim in \textbf{EC}(ii) imply \textbf{LEC}(ii). The second claim in \textbf{EC}(ii) implies \textbf{LEC}(iii).
It is shown in Theorem~\ref{thm:acoe:L}  that the ACOEs hold under
Assumptions~\textbf{W*},  \textbf{B}, and \textbf{LEC}.

\vspace{0.1in}
\noindent
\textbf{Assumption LEC.} There exists a sequence $\{ \a_n \}_{n\in\N^*}$ of nonnegative discount factors such that $\a_n\uparrow1$ as $n\to\infty$ and the following three conditions hold:

(i) the sequence of functions $\{u_{\a_n}\}_{n\in\N^*}$ is lower semi-equicontinuous;

(ii) $\lim_{n\to\infty} u_{\a_n} (x)$ exists for each $x\in\X;$ 

(iii) for each $x\in\X$ and $a\in\A$  the sequence $\{u_{\a_n}\}_{n\in\N^*}$
		is a.u.i. w.r.t. $q(\cdot|x,a).$

\begin{theorem}\label{thm:acoe:L}
Let Assumptions~\textbf{W*} and \textbf{B} hold. Consider a sequence $\{ \a_n\uparrow1 \}_{n\in\N^*}$ of nonnegative discount factors.
If Assumption \textbf{LEC} is satisfied for the sequence $\{ \a_n\}_{n\in\N^*},$
then there exists a stationary  policy $\phi$ such that the ACOEs \eqref{eqn:ACOE} hold with the function $u(x)$
defined in \eqref{eqn:tu:setwise}.
\end{theorem}

\proof{Proof.}
Since Assumptions \textbf{W*} and \textbf{B} hold and $\{ u_{\a_n} \}_{n\in\N^*}$ is lower semi-equicontinuous, then Corollary~\ref{c6.4EF} implies that there exists a stationary policy $\phi$ satisfying \eqref{eqn:ACOI} with $u$ defined in \eqref{eqn:tu:setwise} 
\begin{align}
	w^* + u(x) \geq c(x,\phi(x)) + \int_{\X} u(y)q(dy|x,\phi(x)) .
	\label{eqn:geq}
\end{align}

To prove the ACOEs, it remains to prove the opposite inequality to \eqref{eqn:geq}.
According to Feinberg et al.~\cite[Theorem 2(iv)]{FKZ12}, for each $n\in\N^*$ and $x\in\X$ the discounted-cost optimality equation is
$v_{\a_n}(x) = \min_{a\in \A} [ c(x,a) + \a_n \int_{\X} v_{\a_n} (y) q(dy|x,a) ],$  which, by subtracting $m_{\a}$ from both sides and by replacing $\a_n$ with 1, implies that for all $a\in\A$
\begin{align}
	(1-\a_n)m_{\a_n} + u_{\a_n} (x) \leq 
 c(x,a) +  \int_{\X} u_{\a_n} (y) q(dy|x,a),\qquad x\in\X  .
	\label{eqn:transform dcoe}
\end{align}
Let $n\to \infty.$ In view of \eqref{eqsfkz}, Assumptions~\textbf{LEC}(ii, iii), and Fatou's lemma~\cite[p. 211]{Shi96},
\eqref{eqn:transform dcoe} imply that for all $a\in\A$
\begin{align}
	w^* + u(x) \leq  c(x,a) + \int_{\X} u(y)q(dy|x,a), \qquad x\in\X .
	\label{eq:acoi:rev}
\end{align}
We remark that the integral in \eqref{eqn:transform dcoe} converges to the integral in \eqref{eq:acoi:rev} since the sequence $\{u_{\a_n}\}_{n\in\N^*}$ converges pointwise to $u$ and is u.i.; see Theorem~\ref{thm:uiCondEqui}.
Then, \eqref{eq:acoi:rev} implies
\begin{align}
	w^* + u(x) \leq \min_{a\in \A} [ c(x,a) + \int_{\X}                                                                                                       u(y)q(dy|x,a)]\leq c(x,\phi(x)) + \int_{\X} u (y) q(dy|x,\phi(x)),\quad x\in\X .
	\label{eqn:leq}
\end{align}
 Thus, \eqref{eqn:geq} and \eqref{eqn:leq} imply
\eqref{eqn:ACOE}.
\hfill\Halmos\endproof

In the following example, Assumptions~\textbf{W*}, \textbf{B}, and \textbf{LEC} hold. Hence the ACOEs hold. However, Assumption~\textbf{EC} does not hold.
Therefore, Assumption~\textbf{LEC} is more general than Assumption~\textbf{EC}.

\begin{example}\label{ex:mdpLEC}
{\rm
	Consider $\X = [0,1]$ equipped with the Euclidean metric
	and $\A = \{ a^{(1)} \} .$ The transition probabilities are
	$q(0|x,a^{(1)}) = 1$ for all $x\in\X.$ The cost function is $c(x,a^{(1)}) = \h \{ x \neq 0 \},$ $x\in\X.$ Then the discounted-cost value is $v_\a (x) = u_\a (x) =\h \{ x \neq 0 \},$  $\a\in[0,1)$ and $x\in\X,$ and the average-cost value is $w^* = w_1 (x) = 0,$ $x\in\X.$
	It is straightforward to see that Assumptions~\textbf{W*} and \textbf{B} hold.
	In addition, since the function $u(x) = \h \{ x \neq 0 \}$ is lower semi-continuous, but it is not continuous,
	the sequence of functions $\{u_{\a_n}\}_{n\in\N^*}$ is lower semi-equicontinuous,  but it is not
	equicontinuous for each sequence $\{ \a_n\uparrow1 \}_{n\in\N^*}.$
	Therefore, Assumption~\textbf{LEC} holds since $0\leq u_{\a_n} (x)\leq 1,$ $x\in\X,$
	and Assumption~\textbf{EC} does not hold.
	The \eqref{eqn:ACOE} holds with $w^* = 0,$ $u (x) = \h \{ x \neq 0 \},$ and
	$\phi (x) = a^{(1)},$ $x\in\X.$
	\hfill\Halmos\endproof
}
\end{example}

The following theorem states the validity of ACOEs under
Assumptions~\textbf{S*}, \textbf{B}, and \textbf{LEC}{\rm(ii,iii)}.

\begin{theorem}\label{thm:acoe:S}
Let Assumptions~\textbf{S*} and \textbf{B} hold. Consider a sequence $\{ \a_n\uparrow1 \}_{n\in\N^*}$ of nonnegative discount factors.
If Assumptions~\textbf{LEC}{\rm(ii,iii)} are satisfied for the sequence $\{ \a_n\}_{n\in\N^*},$
then there exists a stationary  policy $\phi$ such that \eqref{eqn:ACOE} holds with the function $u(x)$
defined in \eqref{eqn:tu:setwise}.
\end{theorem}

\proof{Proof.}
According to Theorem~\ref{thm:ACOI:s},
if Assumptions \textbf{S*} and \textbf{B} hold, then we have that: (i) equalities in \eqref{eqsfkz} hold;
(ii) there exists a stationary policy $\phi$ satisfying  ACOI \eqref{eqn:geq} with the function $ u$ defined
in \eqref{eqn:tu:setwise}; and
(iii) for each $n\in\N^*$ and $x\in\X$ the discounted-cost optimality equation is
$v_{\a_n}(x) = \min_{a\in \A} [ c(x,a) + \a_n \int_{\X} v_{\a_n} (y) q(dy|x,a) ].$
Therefore, the same arguments as in the proof of Theorem~\ref{thm:acoe:L} starting from \eqref{eqn:transform dcoe}
imply the validity of \eqref{eqn:ACOE} with $u$ defined
in \eqref{eqn:tu:setwise}.
\hfill\Halmos\endproof

Observe that the MDP described in Example~\ref{ex:mdpLEC}
also satisfies Assumptions~\textbf{S*}, \textbf{B}, and \textbf{LEC}{\rm(ii,iii)}.
We provide Example~\ref{ex:mdpLEC:s},
in which Assumptions~\textbf{S*}, \textbf{B}, and \textbf{LEC}{\rm(ii,iii)} hold.
Hence, the ACOEs hold. However, Assumptions~\textbf{W*}, \textbf{LEC}(i), and \textbf{EC} do not hold.

\begin{example}\label{ex:mdpLEC:s}
{\rm
	Let $\X = [0,1]$ and $\A = \{ a^{(1)} \} .$ The transition probabilities are
	$q(0|x,a^{(1)}) = 1$ for all $x\in\X.$ The cost function is $c(x,a^{(1)}) = D(x),$ where $D$ is the
	Dirichlet function defined as
	\begin{align*}
		D (x) =
		\begin{cases}
			0, & \text{if } x \text{ is rational,} \\
			1, & \text{if } x \text{ is irrational,}
		\end{cases}
		\qquad x\in\X.
	\end{align*}
	Since there is only one available action, Assumption~\textbf{S*} holds.
	The discounted-cost value is $v_\a (x) = u_\a (x) = D(x)=u(x),$  $\a\in[0,1)$ and $x\in\X,$
	and the average-cost value is $w^* = w_1 (x) = 0,$ $x\in\X.$ Then Assumptions~\textbf{B} and \textbf{LEC}{\rm(ii,iii)} hold.
	Hence, the ACOEs \eqref{eqn:ACOE} hold with $w^* = 0,$ $                                                                                                                      u (x) = D(x),$ and $\phi (x) = a^{(1)},$ $x\in\X.$  Thus, the average-cost relative function $u$ is not lower semi-continuous.
	However, since the function $c(x,a^{(1)})=D(x)$ is not lower semi-continuous, Assumption~\textbf{W*} does not hold.
	 Since the function $u(x) =u_\a (x) = D(x)$ is not lower semi-continuous,  Assumptions~\textbf{LEC}(i) and \textbf{EC} do not hold either.
	\hfill\Halmos\endproof
}
\end{example}

\noindent {\bf Acknowledgement.} Research of the first and the third authors was partially supported by NSF grant CMMI-1636193. The authors thank Huizhen (Janey) Yu for valuable remarks.

\renewcommand\refname{\center{\normalfont{\small{REFERENCES}}}}

\end{document}